# Adaptive Control By Regulation-Triggered Batch Least-Squares Estimation of Non-Observable Parameters

**Iasson Karafyllis[*], Maria Kontorinaki[**] and Miroslav Krstic[***]**

[*]Dept. of Mathematics, National Technical University of Athens,
Zografou Campus, 15780, Athens, Greece, email: iasonkar@central.ntua.gr

[**]Department of Statistics and Operations Research, University of Malta,
Tal-Qroqq Campus, Msida, 2080, Malta, email: kontorinmaria@gmail.com

[***]Dept. of Mechanical and Aerospace Eng., University of California, San Diego,
La Jolla, CA 92093-0411, U.S.A., email: krstic@ucsd.edu

**Abstract**

The paper extends a recently proposed indirect, certainty-equivalence, event-triggered adaptive control scheme to the case of non-observable parameters. The extension is achieved by using a novel Batch Least-Squares Identifier (BaLSI), which is activated at the times of the events. The BaLSI guarantees the finite-time asymptotic constancy of the parameter estimates and the fact that the trajectories of the closed-loop system follow the trajectories of the nominal closed-loop system ("nominal" in the sense of the asymptotic parameter estimate, not in the sense of the true unknown parameter). Thus, if the nominal feedback guarantees global asymptotic stability and local exponential stability, then unlike conventional adaptive control, the newly proposed event-triggered adaptive scheme guarantees global asymptotic regulation with a uniform exponential convergence rate. The developed adaptive scheme is tested to a well-known control problem: the state regulation of the wing-rock model. Comparisons with other adaptive schemes are provided for this particular problem.

## 1. Introduction

Adaptive control of linear and nonlinear finite-dimensional systems is an important topic of the control literature. Classical and comprehensive references such as [17,23,24,35] are helpful for the understanding of existing approaches to adaptive control of finite-dimensional systems. Many existing approaches have been also extended to (i) parabolic Partial Differential Equations (PDEs) in one spatial dimension (see [37]), and (ii) hyperbolic PDEs in one spatial dimension (see [1,3] and references therein).

Event-triggered control has attracted considerable attention within the control systems community. Indeed, event-triggered control has been applied to difficult control problems that involve sampling, quantized measurements, output-feedback control, distributed networked control and decentralized control; see [2,4,5,7,8,9,10,11,12,25,26,27,39,40,41,42,45,48]. In all cases, the system under event-triggered control becomes a hybrid dynamical system. Event-triggered direct adaptive control schemes have also appeared in the literature during the last two decades. Event-triggered adaptive control has been applied to globally Lipschitz in the literature of neural networks (see [36,43,47,49]). Direct adaptive control approaches for linear systems have been proposed in [28,29,30,31], where the proposed schemes either employ event-triggering or sampled-data techniques. Event-triggered adaptive control schemes for a special class of nonlinear systems where the input is applied with Zero-Order-Hold were studied in [46]. Adaptive control design methodologies with logic-based switching for linear and nonlinear control systems have been developed in [13,14,15,16,33,34,44] (see also the references therein). The proposed direct supervisory adaptive control schemes in [13,14,15,16,33,34] employ multi-model based estimators of the performance of the "current" controller in conjunction with hierarchical hysteresis switching logic (which is the event-triggered element in the design). Therefore this direct approach is based on an estimation error-triggered controller scheduling.

A different certainty-equivalence, regulation-triggered, indirect adaptive control scheme was proposed in [20] under a parameter observability assumption (but without any persistence of excitation assumption). The adaptive controller in [20] employed a dead-beat, least-squares identifier with delays and allowed the constructive derivation of $KL$ regulation estimates that guarantee the same convergence properties as that of the nominal feedback controller with known parameters. The approach was extended in [21] to the case of reaction-diffusion PDEs in one spatial dimension with constant coefficients.

In the present work, we consider nonlinear systems of the form

$$\dot{x} = f(x,u) + g(x,u)\theta$$
$$x \in \Re^n, u \in \Re^m, \theta \in \Theta \subseteq \Re^l \tag{1.1}$$

where $f : \Re^n \times \Re^m \to \Re^n$, $g : \Re^n \times \Re^m \to \Re^{n \times l}$ are smooth mappings with $f(0,0) = 0$, $g(0,0) = 0$ and $\theta \in \Theta \subseteq \Re^l$ is a vector of constant but unknown parameters that take values in a closed convex set $\Theta \subseteq \Re^l$. By modifying the identifier used in [20], we obtain a new identifier that can work even without any parameter observability assumption. The proposed identifier is a *Batch Least-Squares Identifier* (BaLSI) which is activated at the times of the events and guarantees that:

(i) The parameter estimates $\hat{\theta}(t) \in \Theta \subseteq \Re^l$ change at most $l$ times, where $l$ is the number of the unknown parameters. As a consequence, the parameter estimates remain constant after the time $\tau > 0$ of the last event for which a change in the parameter estimate occurs (*finite-time asymptotic constancy* of the parameter estimates). Moreover, the parameter estimation error $\hat{\theta}(t) - \theta$ satisfies

$$g(x(t),u(t))(\theta - \theta_s) = 0 \text{ for all } t \geq \tau,$$

where $\theta_s$ is the constant value of the parameter estimate after the last event, i.e., $\hat{\theta}(t) = \theta_s$ for $t \geq \tau$.

Moreover, the BaLSI, when combined with a certainty-equivalence controller (as in [20]) achieves that:

(ii) When no change in the parameter estimate occurs, the closed-loop system follows the trajectories of the nominal closed-loop system ("nominal" in the sense of the asymptotic parameter estimate, not in the sense of the true unknown parameter).

To see this, notice that a certainty-equivalence controller $u = k(\hat{\theta}, x)$ applied to system (1.1) with $u = k(\theta, x)$ being the nominal feedback, gives for $t \geq \tau$:

$$\dot{x} = f(x, k(\theta_s, x)) + g(x, k(\theta_s, x))\theta$$
$$= f(x, k(\theta_s, x)) + g(x, k(\theta_s, x))\theta_s + g(x, k(\theta_s, x))(\theta - \theta_s)$$
$$= f(x, k(\theta_s, x)) + g(x, k(\theta_s, x))\theta_s$$

i.e., we follow the trajectories of the nominal closed-loop system $\dot{x} = f(x, k(\theta, x)) + g(x, k(\theta, x))\theta$ with $\theta$ being replaced by $\theta_s$. Therefore, Fact (i) in conjunction with Fact (ii) guarantee that the solution of the closed-loop system presents *exactly the same convergence properties as the nominal closed-loop system*. Thus, if the nominal feedback guarantees global asymptotic stability and local exponential stability, then the proposed event-triggered adaptive scheme guarantees global asymptotic regulation with a uniform exponential convergence rate.

The use of the regulation-triggered schedule of events (as in [20]) allows the following facts:
- (iii) No finite-escape time occurs, even if the nonlinearity is arbitrary.
- (iv) Useful bounds for the solution are obtained, which allow the derivation of $KL$ regulation estimates.
- (v) When no change in the parameter estimate occurs then two consecutive events differ by a constant user-specified time.

Fact (i) in conjunction with Fact (v) guarantee that no Zeno behavior is possible. Finally, the BaLSI guarantees that for many cases, the parameter estimates will converge to the actual values of the parameters (except for a possible set of initial conditions of Lebesgue measure zero). To our knowledge, this collection of desirable features is not exhibited simultaneously by any other adaptive scheme. More specifically, the absence of $KL$ regulation estimates for the supervisory adaptive control schemes in [13,14,15,16,33,34] is explained by the use of the estimation error-triggered policy (instead of our regulation-triggered policy) and the fact that the settling time of the parameter estimate cannot be estimated. This is also true for our scheme (i.e., the time $\tau > 0$ of the last event for which a change in the parameter estimate occurs cannot be estimated) but due to fact (ii) (which does not only hold for the last event but for all events) we are in a position to bound the solution of the closed-loop system by means of an appropriate $KL$ regulation estimate. However, it should be noticed that important robustness results with respect to various errors are provided in [13,14,15,16,33,34], while here we do not consider the possible effect of noise, disturbances and unmodeled dynamics (with the exception of the numerical example in Section V).

In this way, we extend the results contained in [20] to linear and nonlinear finite-dimensional systems with non-observable parameters. The present paper also generalizes the results contained in [20] to systems with parameters that take values in a closed, convex set of the parameter space and consequently the scheme can work even with non-zero parameters (e.g., high-frequency gains).

The structure of the paper is as follows. Section II is devoted to the formulation of the problem and the presentation of the assumptions under which the adaptive regulator is constructed. Section III provides the detailed description of the event-triggered identifier and the adaptive controller. The main results of the present work are given in Section IV (Theorem 4.1, Theorem 4.2, Theorem 4.3 and Corollary 4.4). Section V contains the numerical study of an important illustrative example: the wing-rock model. Section VI contains the proofs of all main results. Finally, the concluding remarks are provided in Section VII.

**Notation.** Throughout this paper, we adopt the following notation.
* For a vector $x \in \Re^n$ we denote by $|x|$ its usual Euclidean norm, by $x'$ its transpose. For a real matrix $A \in \Re^{n \times m}$, $A' \in \Re^{m \times n}$ denotes its transpose and $|A| := \sup\{|Ax| ; x \in \Re^n, |x| = 1\}$ is its induced norm. For a square matrix $A \in \Re^{n \times n}$, $\det(A)$ denotes its determinant and $N(A)$ denotes the null space of $A$, i.e., $N(A) = \{x \in \Re^n : Ax = 0\}$. For a subspace $S$ of $\Re^n$ we denote by $\dim(S)$ its dimension.
* $\Re_+$ denotes the set of non-negative real numbers. $Z_+$ denotes the set of non-negative integers.
* We say that a function $V : \Re^n \to \Re_+$ is positive definite if $V(x) > 0$ for all $x \neq 0$ and $V(0) = 0$. We say that a continuous function $V : \Re^n \to \Re_+$ is radially unbounded if the following property holds: "for every $M > 0$ the set $\{x \in \Re^n : V(x) \leq M\}$ is compact".

* Let $A \subseteq \Re^n$ be an open set, let $U \subseteq \Re^n$ be a set with $A \subseteq U \subseteq cl(A)$, where $cl(A)$ denotes the closure of $A$, and let $\Omega \subseteq \Re$ be a set. By $C^0(U;\Omega)$ we denote the class of continuous mappings on $U$ which take values in $\Omega$. By $C^k(U;\Omega)$, where $k \geq 1$, we denote the class of continuous functions on $U$, which have continuous derivatives of order $k$ on $U$ and also take values in $\Omega$.
* By $K$ we denote the class of strictly increasing $C^0$ functions $a: \Re_+ \to \Re_+$ with $a(0)=0$. By $K_\infty$ we denote the class of strictly increasing $C^0$ functions $a: \Re_+ \to \Re_+$ with $a(0)=0$ and $\lim_{s \to +\infty} a(s) = +\infty$. By $KL$ we denote the set of all continuous functions $\sigma: \Re_+ \times \Re_+ \to \Re_+$ with the properties: (i) for each $t \geq 0$ the mapping $\sigma(\cdot,t)$ is of class $K$; (ii) for each $s \geq 0$, the mapping $\sigma(s,\cdot)$ is non-increasing with $\lim_{t \to +\infty} \sigma(s,t) = 0$.

All stability notions used in this paper are the standard stability notions for time-invariant systems (see [22]).

## 2. Problem Formulation and Assumptions

Consider system (1.1) and suppose that there exist a smooth mapping $k: \Theta \times \Re^n \to \Re^m$ with $k(\theta,0)=0$ for all $\theta \in \Theta$, a constant $\sigma > 0$ and a family of continuously differentiable, positive definite and radially unbounded functions $V(\theta,\cdot) \in C^1(\Re^n;\Re_+)$ parameterized by $\theta \in \Theta$ with the mapping $\Theta \times \Re^n \ni (\theta,x) \to V(\theta,x)$ being continuous, such that the following assumptions hold.

**(H1)** *For each $\theta \in \Theta$, $0 \in \Re^n$ is Globally Asymptotically Stable (GAS) for the closed-loop system*

$$\dot{x} = f(x,k(\theta,x)) + g(x,k(\theta,x))\theta \tag{2.2}$$

*More specifically, the following inequality holds:*

$$\nabla V(\theta,x)\bigl(f(x,k(\theta,x)) + g(x,k(\theta,x))\theta\bigr) \leq -2\sigma V(\theta,x), \text{ for all } \theta \in \Theta, \ x \in \Re^n \tag{2.3}$$

**(H2)** *For every non-empty, compact set $\bar{\Theta} \subseteq \Theta$, the following property holds: "for every $M \geq 0$ there exists $R > 0$ such that the implication $V(\theta,x) \leq M, \theta \in \bar{\Theta} \Rightarrow |x| \leq R$ holds".*

Assumption (H1) is a standard stabilizability assumption (necessary for all possible adaptive control design methodologies). For nonlinear systems, the design of a globally stabilizing state feedback law $u = k(\theta,x)$ is usually performed with the use of a Control Lyapunov Function (CLF, see [6,19,24,38] and references therein). Therefore, the knowledge of the functions $k: \Theta \times \Re^n \to \Re^m$ and $V: \Theta \times \Re^n \to \Re$ is not a demanding requirement. Assumption (H2) is a technical assumption, which requires a "uniform" coercivity property for $V(\theta,\cdot)$ on compact sets of $\Theta \subseteq \Re^l$. Assumption (H2) holds automatically for arbitrary closed, convex sets $\Theta \subseteq \Re^l$ and for functions of the form

$$V(\theta,x) = a_1(\theta,x)x_1^2 + a_2(\theta,x)\bigl(x_2 - \varphi_1(\theta,x_1)\bigr)^2$$
$$+ \ldots + a_n(\theta,x)\bigl(x_n - \varphi_{n-1}(\theta,x_1,\ldots,x_{n-1})\bigr)^2$$

where $a_i: \Theta \times \Re^n \to [\kappa,+\infty)$ ($i=1,\ldots,n$) are functions bounded from below by a positive constant $\kappa > 0$ and $\varphi_i: \Theta \times \Re^i \to \Re$ ($i=1,\ldots,n-1$) are continuous functions with $\varphi_i(\theta,0)=0$ for all $\theta \in \Theta$ and $i=1,\ldots,n-1$. The above functional form is met frequently in the study of nonlinear triangular single-input systems of the form $\dot{x}_1 = f_1(\theta,x_1,x_2)$, ...., $\dot{x}_{n-1} = f_{n-1}(\theta,x_1,\ldots,x_n)$, $\dot{x}_n = f_n(\theta,x_1,\ldots,x_n,u)$.

In [20] we used the following parameter observability assumption:

**(H3)** *There exists a positive integer $N$ such that the following implication holds:*

"*For every set of $N$ times $0 = \tau_0 < \tau_1 < ... < \tau_N$ and for every vectors $\theta, d_0, ..., d_N \in \Theta$ with $d_i \neq 0$ for $i = 0,...,N$ the only right differentiable mapping $x \in C^0([0,\tau_N]; \Re^n) \cap C^1([0,\tau_N] \setminus \{\tau_0,...,\tau_N\}; \Re^n)$ satisfying $\dot{x}(t) = f(x(t), k(\theta + d_i, x(t))) + g(x(t), k(\theta + d_i, x(t)))\theta$ for $t \in [\tau_i, \tau_{i+1})$, $i = 0,...,N-1$, $g(x(t), k(\theta + d_j, x(t)))d_{i+1} = 0$ for all $t \in [\tau_j, \tau_{j+1}]$, $i = 0,...,N-1$, $j = 0,...,i$, is the identically zero mapping, i.e., $x(t) = 0$ for all $t \in [0, \tau_N]$.*"

Assumption (H3) was used in [20] in order to guarantee finite-time identification of the parameters. However, in what follows *we will not employ* Assumption (H3) and consequently we won't be able to guarantee finite-time identification. Neither Assumption (H3) nor finite-time identification are encountered in conventional adaptive control.

## 3. Event-Triggered Identifier for a Certainty-Equivalence Adaptive Controller

In this section we introduce the adaptive control law. The reader interested in a quick access to the adaptive controller may immediately refer to (3.1), (3.2), (3.4), (3.5), (3.13) and then resume reading the rest of this section for explanations.

The control action between two consecutive events is governed by the nominal feedback $u = k(\theta, x)$ with the unknown $\theta \in \Theta$ replaced by its estimate $\hat{\theta}$. Moreover, the estimate $\hat{\theta}$ of the unknown $\theta \in \Theta$ is kept constant between two consecutive events. In other words, we have

$$\begin{aligned} u(t) &= k(\hat{\theta}(\tau_i), x(t)) \quad,\quad t \in [\tau_i, \tau_{i+1}), i \in Z_+ \\ \hat{\theta}(t) &= \hat{\theta}(\tau_i) \quad,\quad t \in [\tau_i, \tau_{i+1}), i \in Z_+ \end{aligned} \qquad (3.1)$$

where $\{\tau_i \geq 0\}_{i=0}^{\infty}$ is the sequence of times of the events that satisfies

$$\begin{aligned} \tau_{i+1} &= \min(\tau_i + T, r_i) \quad,\quad i \in Z_+ \\ \tau_0 &= 0 \end{aligned} \qquad (3.2)$$

where $T > 0$ is a positive constant (one of the tunable parameters of the proposed scheme) and $r_i > \tau_i$ is a time instant determined by the event trigger.

Let $a \in C^0(\Re^n; \Re_+)$ be a positive definite function (the second tunable parameter of the proposed scheme). The event trigger sets $r_i > \tau_i$ to be the smallest time $t > \tau_i$ for which

$$V(\hat{\theta}(\tau_i), x(t)) = V(\hat{\theta}(\tau_i), x(\tau_i)) + a(x(\tau_i)) \qquad (3.3)$$

where $x(t)$ denotes the solution of (1.1) with $u(t) = k(\hat{\theta}(\tau_i), x(t))$. For the case that a time $t > \tau_i$ satisfying (3.3) does not exist, we set $r_i = +\infty$. For the case $x(\tau_i) = 0$ we set $r_i := \tau_i + T$. Formally, the event trigger is described by the equations:

$$r_i := \inf\left\{ t > \tau_i : V(\hat{\theta}(\tau_i), x(t)) = V(\hat{\theta}(\tau_i), x(\tau_i)) + a(x(\tau_i)) \right\}, \text{ for } x(\tau_i) \neq 0 \qquad (3.4)$$

$$r_i := \tau_i + T, \text{ for } x(\tau_i) = 0 \qquad (3.5)$$

The description of the adaptive control scheme is completed by the parameter update law, which is activated at the times of the events.

In order to estimate the unknown vector $\theta \in \Theta$, we develop the Batch Least-Squares Identifier (BaLSI). Notice that (by virtue of (1.1)) for every $t, s \geq 0$ the following equation holds:

$$x(t) - x(s) = \int_s^t f(x(r), u(r)) dr + \left( \int_s^t g(x(r), u(r)) dr \right) \theta \quad (3.6)$$

Define for every $i \in Z_+$ the function $h_i : \Re^l \to \Re_+$ by the formula

$$h_i(\vartheta) := \int_0^{\tau_{i+1}} \int_0^{\tau_{i+1}} |p(t,s) - q(t,s)\vartheta|^2 \, ds \, dt \quad (3.7)$$

where

$$p(t,s) := x(t) - x(s) - \int_s^t f(x(r), u(r)) dr \quad (3.8)$$

$$q(t,s) := \int_s^t g(x(r), u(r)) dr \quad (3.9)$$

It follows from (3.6) and (3.7) that for every $i \in Z_+$ the function $h_i(\vartheta)$ has a global minimum at $\vartheta = \theta$ with $h_i(\theta) = 0$. Consequently, we get from Fermat's theorem for extrema that the following equation holds:

$$\int_0^{\tau_{i+1}} \int_0^{\tau_{i+1}} q'(t,s) p(t,s) ds \, dt = \left( \int_0^{\tau_{i+1}} \int_0^{\tau_{i+1}} q'(t,s) q(t,s) ds \, dt \right) \theta \quad (3.10)$$

It should be noticed that the matrix $G(\tau) := \left( \int_0^\tau \int_0^\tau q'(t,s) q(t,s) ds \, dt \right) \in \Re^{l \times l}$ is symmetric and positive semi-definite. Consequently, if $G(\tau_{i+1}) \in \Re^{l \times l}$ is invertible (i.e., $\det(G(\tau_{i+1})) \neq 0$) then $G(\tau_{i+1}) \in \Re^{l \times l}$ is positive definite with $\det(G(\tau_{i+1})) > 0$ and

$$\theta = \left( \int_0^{\tau_{i+1}} \int_0^{\tau_{i+1}} q'(t,s) q(t,s) ds \, dt \right)^{-1} \int_0^{\tau_{i+1}} \int_0^{\tau_{i+1}} q'(t,s) p(t,s) ds \, dt \quad (3.11)$$

The estimate (3.11) is nothing else but the least squares estimate of the unknown vector $\theta \in \Re^l$ on the interval $[0, \tau_{i+1}]$. In the general case, the following convex optimization problem with linear equality constraints

$$\min_{\vartheta \in \Theta} |\vartheta - \hat{\theta}(\tau_i)|^2$$

$$\text{s.t.} \quad (3.12)$$

$$\int_0^{\tau_{i+1}} \int_0^{\tau_{i+1}} q'(t,s) p(t,s) ds \, dt = \left( \int_0^{\tau_{i+1}} \int_0^{\tau_{i+1}} q'(t,s) q(t,s) ds \, dt \right) \vartheta$$

has a unique solution. We can therefore define the following parameter update law:

$$\hat{\theta}(\tau_{i+1}) = \arg\min \left\{ |\vartheta - \hat{\theta}(\tau_i)|^2 : \vartheta \in \Theta, \int_0^{\tau_{i+1}} \int_0^{\tau_{i+1}} q'(t,s) p(t,s) ds \, dt = \left( \int_0^{\tau_{i+1}} \int_0^{\tau_{i+1}} q'(t,s) q(t,s) ds \, dt \right) \vartheta \right\} \quad (3.13)$$

which is the BaLSI. It should also be noticed that the operator involved in (3.13) is not a continuous operator. However, in practice an accurate continuous approximation of the parameter update law (3.13) may be used.

**Remark 3.1: (a)** The parameter update law (3.13) (the BaLSI) is the key difference of the proposed scheme and the scheme in [20]. More specifically, in [20] the least-squares identifier used a parameter update law of the form (3.13) for which the lower limits of the integrals appearing in (3.13) were not necessarily zero.

**(b)** The BaLSI can be implemented by a set of ODEs. Indeed, an implementation of the parameter update law (3.13) is given by the following ODEs

$$\begin{aligned} \dot{z} &= Cx \quad, \quad z \in \Re^j \\ \dot{B} &= t\bigl(Cg(x,u)\bigr)' \quad, \quad B \in \Re^{l \times j} \\ \dot{w} &= \bigl(Cg(x,u)\bigr)'(z+\phi) + BCf(x,u) \quad, \quad w \in \Re^l \\ \dot{\phi} &= tCf(x,u) \quad, \quad \phi \in \Re^j \\ \dot{Y} &= 2(BCx - w) \quad, \quad Y \in \Re^l \\ \ddot{Q} &= 2\bigl(BCg(x,u)\bigr)' + 2BCg(x,u) \quad, \quad Q \in \Re^{l \times l} \end{aligned} \quad (3.14)$$

where $C \in \Re^{j \times n}$ is a constant matrix with $1 \le j \le n$, $rank(C) = j$ such that $\bigl(I - C'(CC')^{-1}C\bigr)g(x,u) = 0$ for all $(x,u) \in \Re^{n \times m}$, with initial conditions $z(0) = \phi(0) = 0$, $Q(0) = \dot{Q}(0) = 0$, $B(0) = 0$, $Y(0) = w(0) = 0$. The parameter update law (3.13) is given by:

$$\hat{\theta}(\tau_{i+1}) = \arg\min\left\{ \bigl|\vartheta - \hat{\theta}(\tau_i)\bigr|^2 : \vartheta \in \Theta, Y(\tau_{i+1}) = Q(\tau_{i+1})\vartheta \right\} \quad (3.15)$$

and notice that

$$\begin{aligned} Y(t) &= \int_0^t\int_0^t q'(r,s)C'Cp(r,s)\,ds\,dr \\ Q(t) &= \int_0^t\int_0^t q'(r,s)C'Cq(r,s)\,ds\,dr \\ B(t) &= \int_0^t q'(t,s)C'\,ds \\ w(t) &= \int_0^t q'(t,s)C'C\left(x(s) + \int_s^t f(x(r),u(r))\,dr\right)ds \\ \phi(t) &= \int_0^t\left(\int_s^t Cf(x(r),u(r))\,dr\right)ds \end{aligned}$$

Moreover, using the fact that $Q \in \Re^{l \times l}$ is symmetric, it is possible to use only $l(l+1)/2$ from the $l^2$ second order ODEs $\ddot{Q} = 2\bigl(Cg(x,u)\bigr)'B' + 2BCg(x,u)$. However, in many cases, the structure of the control system (1.1) allows a large reduction of the number of ODEs that are needed for the implementation of the parameter update law given by (3.13), by selecting the matrix $C \in \Re^{j \times n}$ in an appropriate way.

It is straightforward to show (using (3.14)) that if local exponential regulation of $x$ is achieved then the variables $z \in \Re^j, B \in \Re^{l \times j}, w \in \Re^l, \phi \in \Re^j$ remain bounded for all $t \ge 0$. Moreover, in this case it may be shown that $|Q(t)| = O(t)$ and $|Y(t)| = O(t)$, i.e., the entries of $Y, Q$ are bounded by a linear time function. For practical operation, the system (3.14) may be re-initiated frequently in order to keep the entries of the matrix $Q \in \Re^{l \times l}$ small (and consequently the components of the vector $Y \in \Re^l$, since we always have $Y(t) = Q(t)\theta)$ when $Q \in \Re^{l \times l}$ is non-singular: in this case the exact value of $\theta \in \Re^l$ has been found and the updates are used only in case that a change of the parameter values occurs.

In practice, it is better to avoid the implementation of the parameter update law (3.13), because due to the presence of modeling and measurement errors the equation $Y(\tau_{i+1}) = Q(\tau_{i+1})\theta$ (guaranteed by (3.10) when modeling and measurement errors are absent) may not hold. Therefore, there is no guarantee that the set $\{\vartheta \in \Theta : Y(\tau_{i+1}) = Q(\tau_{i+1})\vartheta\}$ is non-empty. Consequently, we may need to relax the minimization problem (3.15) and use the following parameter update law instead of (3.15):

$$\hat{\theta}(\tau_{i+1}) = \arg\min\left\{ \left|\vartheta - \hat{\theta}(\tau_i)\right|^2 + \beta \left|Y(\tau_{i+1}) - Q(\tau_{i+1})\vartheta\right|^2 : \vartheta \in \Theta \right\} \qquad (3.16)$$

where $\beta > 0$ is a large positive constant. The parameter update law (3.16) has additional advantages compared to (3.15), since (3.16) introduces a regularization effect. To see this, notice that when $\Theta = \Re^l$, (3.16) gives $\hat{\theta}(\tau_{i+1}) = \left(\beta^{-1}I + Q^2(\tau_{i+1})\right)^{-1}\left(\beta^{-1}\hat{\theta}(\tau_i) + Q(\tau_{i+1})Y(\tau_{i+1})\right)$ where (due to the fact that $Q(\tau_{i+1})$ is a symmetric and positive semi-definite matrix) the matrix $\beta^{-1}I + Q^2(\tau_{i+1})$ is always positive definite.

## 4. Statements of Stability Results

We consider the plant (1.1) with the controller (3.1), (3.2), (3.4), (3.5) and the parameter estimator (3.13). The first main result guarantees global regulation of $x$ to zero.

**Theorem 4.1:** *Consider the control system (1.1) under assumptions (H1), (H2). Let $T > 0$ be a positive constant and let $a: \Re^n \to \Re_+$ be a continuous, positive definite function. Then there exists a family of KL mappings $\omega_{\theta,\hat{\theta}} \in KL$ parameterized by $\theta \in \Theta$, $\hat{\theta} \in \Theta$ such that for every $\theta \in \Theta$, $x_0 \in \Re^n$, $\hat{\theta}_0 \in \Theta$ the solution of the hybrid closed-loop system (1.1) with (3.1), (3.2), (3.4), (3.5), (3.13) and initial conditions $x(0) = x_0$, $\hat{\theta}(0) = \hat{\theta}_0$ is unique, is defined for all $t \geq 0$ and satisfies $|x(t)| \leq \omega_{\theta,\hat{\theta}}(|x_0|,t)$ for all $t \geq 0$. Moreover, there exist $\tau \geq 0$, $\theta_s \in \Theta$ (both depending on $\theta \in \Theta$, $x_0 \in \Re^n$, $\hat{\theta}_0 \in \Theta$) such that $\hat{\theta}(t) = \theta_s$ for all $t \geq \tau$ and the equation $g(x(t),u(t))(\theta - \theta_s) = 0$ holds for all $t \geq 0$. Finally, if assumption (H3) holds and $x_0 \neq 0$ then $\hat{\theta}(t) = \theta$ for all $t \geq NT$.*

**Remarks on Theorem 4.1:** **(a)** Theorem 4.1 guarantees that there is a finite settling time $\tau \geq 0$ for the parameter estimate $\hat{\theta}(t) \in \Theta$. Unfortunately, an upper bound of the settling time cannot be provided. **(b)** The condition $g(x(t),u(t))(\theta - \theta_s) = 0$ for all $t \geq 0$, for many cases can only be satisfied on a (Lebesgue) measure zero set of initial conditions $x_0 \in \Re^n$, $\hat{\theta}_0 \in \Theta$. Therefore, in such cases, identification of the parameters is not possible only for a measure zero set of initial conditions $x_0 \in \Re^n$, $\hat{\theta}_0 \in \Theta$. **(c)** It is important to notice that no assumption for persistency of excitation is made in Theorem 4.1. **(d)** The proof of Theorem 4.1 shows that at most $l$ switchings of the value of the parameter estimate $\hat{\theta}(t)$ can occur. Moreover, the estimate $\tau_i \geq (i-l)T$ holds for all $i \geq l$, indicating that the times of the events cannot have a finite accumulation point (Zeno behavior).

The second main result guarantees local exponential regulation of $x$ to zero under the assumption that the nominal feedback law $u = k(\theta, x)$ achieves local exponential stabilization.

**Theorem 4.2:** *Consider the control system (1.1) under assumptions (H1), (H2). Moreover, suppose that for each $\theta \in \Theta$, $0 \in \Re^n$ is Locally Exponentially Stable (LES) for the closed-loop system (2.2) and that for every nonempty, compact set $\bar{\Theta} \subseteq \Theta$ there exist constants $R > 0$, $K_2 > K_1 > 0$ such that*

$$K_1 |x|^2 \leq V(\theta, x) \leq K_2 |x|^2, \text{ for all } x \in \Re^n, \theta \in \bar{\Theta} \text{ with } |x| \leq R \quad (4.1)$$

*Let $T > 0$ be a positive constant and let $a: \Re^n \to \Re_+$ be a continuous, positive definite function that satisfies $\sup\{|x|^{-2} a(x): x \in \Re^n, x \neq 0, |x| \leq \delta\} < +\infty$ for certain $\delta > 0$. Then there exists a family of constants $M_{\theta, \hat{\theta}}, \bar{R}_{\theta, \hat{\theta}} > 0$ parameterized by $(\theta, \hat{\theta}) \in \Theta \times \Theta$, such that for every $\theta \in \Theta$, $x_0 \in \Re^n$, $\hat{\theta}_0 \in \Theta$ with $|x_0| \leq \bar{R}_{\theta, \hat{\theta}_0}$ the solution of the hybrid closed-loop system (1.1) with (3.1), (3.2), (3.4), (3.5), (3.13) and initial conditions $x(0) = x_0$, $\hat{\theta}(0) = \hat{\theta}_0$ satisfies the estimate $|x(t)| \leq M_{\theta, \hat{\theta}_0} \exp(-\sigma t)|x_0|$ for all $t \geq 0$, with $\sigma > 0$ being the constant involved in (2.3).*

It should be noticed that Theorem 4.2 guarantees that the local exponential stability estimate $|x(t)| \leq M_{\theta, \hat{\theta}_0} \exp(-\sigma t)|x_0|$ holds when $|x_0| \leq \bar{R}_{\theta, \hat{\theta}_0}$ and for arbitrary initial condition $\hat{\theta}_0 \in \Theta$. In other words, the adjective "local" refers only to $x$ and not to $\hat{\theta}$. Moreover, the reader should notice that the event-triggered adaptive scheme (3.1), (3.2), (3.4), (3.5), (3.13) guarantees convergence with the same convergence rate $\sigma > 0$ as the nominal feedback controller with known parameter values. This is not possible for conventional adaptive control: in conventional adaptive control there is no uniform exponential convergence rate for all initial conditions. The uniform exponential convergence rate is achieved by estimating the parameter vector in an appropriate way. To see this notice that when the parameter estimate satisfies $\hat{\theta}(t) = \theta_s$ for all $t \geq \tau$, the hybrid closed-loop system (1.1) with (3.1), (3.2), (3.4), (3.5), (3.13) is described for $t \geq \tau$ by the equation

$$\begin{aligned}\dot{x} &= f(x, k(\theta_s, x)) + g(x, k(\theta_s, x))\theta \\ &= f(x, k(\theta_s, x)) + g(x, k(\theta_s, x))\theta_s + g(x, k(\theta_s, x))(\theta - \theta_s)\end{aligned}$$

The least-squares identifier guarantees that the parameter estimate satisfies $g(x(t), k(\theta_s, x(t)))(\theta - \theta_s) = 0$ for all $t \geq \tau$ and consequently, the hybrid closed-loop system (1.1) with (3.1), (3.2), (3.4), (3.5), (3.13) is described for $t \geq \tau$ by the equation

$$\dot{x} = f(x, k(\theta_s, x)) + g(x, k(\theta_s, x))\theta_s$$

i.e., as if the parameter vector were $\theta_s$ instead of $\theta$. Therefore, the solution of the adaptive closed-loop system coincides with the solution of (2.2) with $\theta$ replaced by $\theta_s$, i.e., the solution of the closed-loop system (1.1) with $u = k(\theta, x)$ and known parameter values.

Finally, it should be emphasized that in addition to the exponential regulation estimate $|x(t)| \leq M_{\theta, \hat{\theta}_0} \exp(-\sigma t)|x_0|$, Theorem 4.2 guarantees all the conclusions of Theorem 4.1 (because all assumptions of Theorem 4.1 are fulfilled).

The following result guarantees global exponential convergence of $x$ to zero under the assumption that the nominal feedback law $u = k(\theta, x)$ achieves global exponential stabilization.

**Theorem 4.3:** *Consider system (1.1) under assumptions (H1), (H2). Moreover, suppose that for each $\theta \in \Theta$, $0 \in \Re^n$ is Globally Exponentially Stable (GES) for the closed-loop system (2.2) and that for every nonempty, compact set $\bar{\Theta} \subseteq \Theta$ there exist constants $K_2 > K_1 > 0$ such that*

$$K_1 |x|^2 \leq V(\theta, x) \leq K_2 |x|^2, \text{ for all } x \in \Re^n, \theta \in \bar{\Theta} \tag{4.2}$$

*Let $T > 0$ be a positive constant and let $a : \Re^n \to \Re_+$ be a continuous, positive definite function that satisfies $\sup\{|x|^{-2} a(x) : x \in \Re^n, x \neq 0\} < +\infty$. Then there exists a family of constants $M_{\theta, \hat{\theta}} > 0$ parameterized by $\theta \in \Theta$, $\hat{\theta} \in \Theta$, such that for every $\theta \in \Theta$, $x_0 \in \Re^n$, $\hat{\theta}_0 \in \Theta$ the solution of the hybrid closed-loop system (1.1) with (3.1), (3.2), (3.4), (3.5), (3.13) and initial conditions $x(0) = x_0$, $\hat{\theta}(0) = \hat{\theta}_0$ satisfies the estimate $|x(t)| \leq M_{\theta, \hat{\theta}_0} \exp(-\sigma t) |x_0|$ for all $t \geq 0$, with $\sigma > 0$ being the constant involved in (2.3).*

Finally, the following corollary deals with the case of controllable Linear Time-Invariant (LTI) single-input systems with unknown parameters.

**Corollary 4.4:** *Consider the system*

$$\begin{aligned} \dot{x}_i &= \sum_{j=1}^{i} \theta_{i,j} x_j + \theta_{i,i+1} x_{i+1}, i = 1, \ldots, n \\ x &= (x_1, \ldots, x_n)' \in \Re^n, u = x_{n+1} \in \Re, \\ \theta &= (\theta_{1,1}, \theta_{1,2}, \ldots, \theta_{n,n+1})' \in \Theta = \left\{ \vartheta \in \Re^l : \vartheta_{i,i+1} \geq \kappa, i = 1, \ldots, n \right\} \\ l &= \frac{n(n+3)}{2} \end{aligned} \tag{4.3}$$

*where $l = \dfrac{n(n+3)}{2}$ and $\kappa > 0$ is a constant. Then for every $\sigma > 0$, $\theta \in \Theta$ there exist a symmetric, positive definite matrix $P(\theta) = \{p_{i,j}(\theta) : i, j = 1, \ldots, n\} \in \Re^{n \times n}$ and a vector $\tilde{k}(\theta) \in \Re^n$ such that the inequality $\sum_{i=1}^{n}\sum_{j=1}^{n} p_{i,j}(\theta) \dot{x}_i x_j + \sum_{i=1}^{n}\sum_{j=1}^{n} p_{i,j}(\theta) x_i \dot{x}_j \leq -2\sigma \sum_{i=1}^{n}\sum_{j=1}^{n} p_{i,j}(\theta) x_i x_j$ with $u = \tilde{k}'(\theta) x$ holds for all $x \in \Re^n$. Moreover, the mappings $\Theta \ni \theta \to P(\theta) \in \Re^{n \times n}$, $\Theta \ni \theta \to \tilde{k}(\theta) \in \Re^n$ are continuous. Finally, for every $T > 0$, $A > 0$, there exists a family of constants $M_{\theta, \hat{\theta}} > 0$ parameterized by $\theta \in \Theta$, $\hat{\theta} \in \Theta$, such that for all $\theta \in \Theta$, $x_0 \in \Re^n$, $\hat{\theta}_0 \in \Theta$ the solution of the hybrid closed-loop system (1.1) with*

$$f(x, u) \equiv 0, \quad g(x, u)\theta := \begin{bmatrix} \theta_{1,1} x_1 + \theta_{1,2} x_2 \\ \vdots \\ \theta_{n,1} x_1 + \ldots + \theta_{n,n} x_n + \theta_{n,n+1} u \end{bmatrix}, \quad a(x) := A|x|^2, \quad V(\theta, x) := x' P(\theta) x, \quad k(\theta, x) := \tilde{k}'(\theta) x,$$

*(3.1), (3.2), (3.4), (3.5), (3.13) and initial conditions $x(0) = x_0$, $\hat{\theta}(0) = \hat{\theta}_0$ satisfies the estimate $|x(t)| \leq M_{\theta, \hat{\theta}_0} \exp(-\sigma t) |x_0|$ for all $t \geq 0$.*

# 5. Adaptive Control of the Wing-Rock Model

The wing-rock model proposed and used in [24,32] with zero torque at equilibrium is given by the system

$$\begin{aligned}
\dot{x}_1 &= x_2 \\
\dot{x}_2 &= \theta_1 x_1 + \theta_2 x_2 + \theta_3 |x_1| x_2 + \theta_4 |x_2| x_2 + \theta_5 x_3 \\
\dot{x}_3 &= -\mu x_3 + \mu u \\
x &= (x_1, x_2, x_3)' \in \Re^3, u \in \Re \\
\theta &= (\theta_1, \theta_2, \theta_3, \theta_4, \theta_5)' \in \Theta = \Re^4 \times [\kappa, +\infty)
\end{aligned} \quad (5.1)$$

where $\kappa > 0$ and $\mu > 0$ are known parameters. A locally Lipschitz nonlinear feedback law that achieves local exponential stabilization and global asymptotic stabilization of the equilibrium point $0 \in \Re^3$ is given by the formula:

$$\begin{aligned}
k(\theta, x) &= x_3 - \mu^{-1}\theta_5 \left(x_2 + Lx_1\right) + \mu^{-1}\frac{\partial \varphi}{\partial x_1}(\theta, x_1, x_2)x_2 \\
&+ \mu^{-1}\frac{\partial \varphi}{\partial x_2}(\theta, x_1, x_2)\left(\theta_1 x_1 + \theta_2 x_2 + \theta_3 |x_1| x_2 + \theta_4 |x_2| x_2 + \theta_5 x_3\right) \\
&- \mu^{-1}L(x_3 - \varphi(\theta, x_1, x_2))
\end{aligned} \quad (5.2)$$

where $L > 1$ is a constant and

$$\begin{aligned}
\varphi(\theta, x_1, x_2) &:= -\theta_5^{-1}\left((1 + \theta_1)x_1 + \theta_2 x_2\right) \\
&\quad - \theta_5^{-1}\left(L + \beta(\theta, x_1, x_2) + \frac{L^2}{4}(\beta(\theta, x_1, x_2))^2\right)(x_2 + Lx_1)
\end{aligned} \quad (5.3)$$

$$\beta(\theta, x_1, x_2) := 1 + \frac{1}{2}\theta_3^2 x_1^2 + \frac{1}{2}\theta_4^2 x_2^2 \quad (5.4)$$

More specifically, the feedback law $u = k(\theta, x)$ guarantees the differential inequality

$$\dot{V}(\theta, x) \leq -2(L-1)V(\theta, x), \text{ for all } \theta \in \Theta, x \in \Re^3 \quad (5.5)$$

for the Control Lyapunov Function

$$V(\theta, x) := \frac{1}{2}x_1^2 + \frac{1}{2}(x_2 + Lx_1)^2 + \frac{1}{2}(x_3 - \varphi(\theta, x_1, x_2))^2 \quad (5.6)$$

Therefore, Assumption (H1) and Assumption (H2) hold with $\sigma := L - 1$. Moreover, for every nonempty, compact set $\bar{\Theta} \subseteq \Theta$ there exist constants $R > 0$, $K_2 > K_1 > 0$ such that (4.1) holds. Therefore all assumptions of Theorem 4.2 hold for system (5.1). Thus, for every $\kappa, \mu, T > 0$, $a \in C^0(\Re^3; \Re_+)$ being a positive definite function that satisfies $\sup\left\{|x|^{-2} a(x) : x \in \Re^3, x \neq 0, |x| \leq \delta\right\} < +\infty$ for certain $\delta > 0$, there exists a family of constants $M_{\theta,\hat{\theta}}, \bar{R}_{\theta,\hat{\theta}} > 0$ parameterized by $(\theta, \hat{\theta}) \in \Theta \times \Theta$, such that for every $\theta \in \Theta$, $x_0 \in \Re^3$, $\hat{\theta}_0 \in \Theta$ with $|x_0| \leq \bar{R}_{\theta,\hat{\theta}_0}$ the solution of the hybrid closed-loop system (5.1) with (3.1), (3.2), (3.4), (3.5), (3.13) and initial conditions $x(0) = x_0$, $\hat{\theta}(0) = \hat{\theta}_0$ satisfies the estimate $|x(t)| \leq M_{\theta,\hat{\theta}_0} \exp(-\sigma t)|x_0|$ for all $t \geq 0$, with $\sigma := L - 1 > 0$ and $\Theta = \Re^4 \times [\kappa, +\infty)$.

We studied numerically system (5.1) with

$$\theta_1 = -26.67, \; \theta_2 = 0.76485, \; \theta_3 = -2.9225, \; \theta_4 = 0, \; \theta_5 = 1.5, \; \mu = 15, \; \kappa = 1$$

which are exactly the same parameter values used in [24]. Notice that $0 \in \Re^3$ is unstable for the open-loop system (5.1) with $u \equiv 0$ and the solution is attracted by a limit cycle (the wing-rock phenomenon; see [24]). The controller parameters were selected as follows:

$$L = 1.5, \; T = 0.4, \; a(x) := 2 \cdot 10^5 \left(|x|^2 + |x|^4\right), \text{ for } x \in \Re^3$$

and following Remark 3.1(b) the hybrid adaptive controller was implemented by using (3.16) with $\beta = 10^{17}$ and the set of ODEs:

$$\dot{z} = x_2 \tag{5.7}$$

$$\dot{Y}_i = 2(b_i x_2 - w_i), \; i = 1,...,5 \tag{5.8}$$

$$\begin{aligned} \dot{b}_i &= t\zeta_i(x) \\ \dot{w}_i &= z\zeta_i(x) \end{aligned}, \; i = 1,...,5 \tag{5.9}$$

$$\ddot{Q}_{i,j} = 2(b_j \zeta_i(x) + b_i \zeta_j(x)), \; j = 1,...,i, \; i = 1,...,5 \tag{5.10}$$

with initial conditions $z(0) = Y_i(0) = b_i(0) = w_i(0) = Q_{i,j}(0) = \dot{Q}_{i,j}(0) = 0$ for $i, j = 1,...,5$, where

$$\begin{aligned} \zeta_1(x) &:= x_1, & \zeta_2(x) &:= x_2, & \zeta_5(x) &:= x_3 \\ \zeta_3(x) &:= |x_1|x_2, & \zeta_4(x) &:= |x_2|x_2 \end{aligned} \tag{5.11}$$

for $x \in \Re^3$.

In order to compare the performance of the closed-loop system (5.1) with (3.1), (3.2), (3.4), (3.5), (3.16), we also used the adaptive controller based on the extended matching design studied for general nonlinear systems in Chapter 3 and Chapter 4 in [24]. More specifically, the extended matching design (combined with the projection schemes explained in Appendix E of the book [24]) gives the following adaptive controller

$$\frac{d}{dt}\hat{\theta}_i = \gamma \frac{\partial V}{\partial x_2}(\hat{\theta}, x)\zeta_i(x), \; i = 1,...,4$$

$$\frac{d}{dt}\hat{\theta}_5 = \gamma \frac{\partial V}{\partial x_2}(\hat{\theta}, x)\zeta_5(x) \begin{cases} 1, & \text{if } \hat{\theta}_5 \geq \kappa \text{ or } \frac{\partial V}{\partial x_2}(\hat{\theta}, x)\zeta_5(x) \geq 0 \\ 1 - \min(1, \varepsilon^{-1}(\kappa - \hat{\theta}_5)), & \text{if otherwise} \end{cases} \tag{5.12}$$

$$u = k(\hat{\theta}, x) + \mu^{-1}\gamma \frac{\partial V}{\partial x_2}(\hat{\theta}, x)\frac{\partial \varphi}{\partial \hat{\theta}}(\hat{\theta}, x_1, x_2)\zeta(x)$$

where $V$ is given by (5.6), $k$ is given by (5.2), (5.3), (5.4) with $L = 1.5$, $\gamma > 0$ is a constant (the adaptation gain),

$$\zeta(x) := (\zeta_1(x), \zeta_2(x), \zeta_3(x), \zeta_4(x), \zeta_5(x))', \text{ for } x \in \Re^3. \tag{5.13}$$

and $\varepsilon \in (0, \kappa)$ is the constant for which the inequality $\hat{\theta}_5(t) \geq \kappa - \varepsilon$ is guaranteed to hold for all $t \geq 0$. The adaptation gain was selected $\gamma = 10$ and $\varepsilon \in (0, \kappa)$ was selected to be equal to 0.001.

We compared the performance for many initial conditions and representative results are presented next. Although, we are not able to verify Assumption (H3) for system (5.1) with nominal controller given by (5.2), finite-time exact estimation of the parameters was achieved for all tested initial conditions in a few (two or three) triggers.

The results shown in Figures 1-9 are obtained by using two different initial conditions for system (5.1):

- 1st initial condition: $x_0 = (-0.35, -0.5, 0.05)'$.
- 2nd initial condition: $x_0 = (0.4, 0, 0)'$.

The 2nd initial condition is exactly the initial condition for which numerical results are presented in [24]. In all cases, the initial condition for the parameter estimates is (as in [24])

$$\hat{\theta}_0 = 1.35\theta$$

Figures 1, 2 show the projected trajectories of the solutions on the $x_1 - x_2$ plane for the 1st and 2nd initial condition, respectively, when no measurement errors are present. It is clearly shown that the trajectories of the solutions of the closed-loop system with the event-triggered adaptive scheme differ significantly from the solutions of the closed-loop system with the extended-matching design. The trajectories of the closed-loop system with the event-triggered adaptive scheme (3.1), (3.2), (3.4), (3.5), (3.16) are similar to the trajectories of the nominal closed-loop system (5.1) with (5.2) due to the fact that the exact values of the parameters are found after an initial transient period. This fact is shown in Figures 7, 8, where it is also shown that the parameter estimates for the extended-matching design fail to converge to the actual values of the parameters.

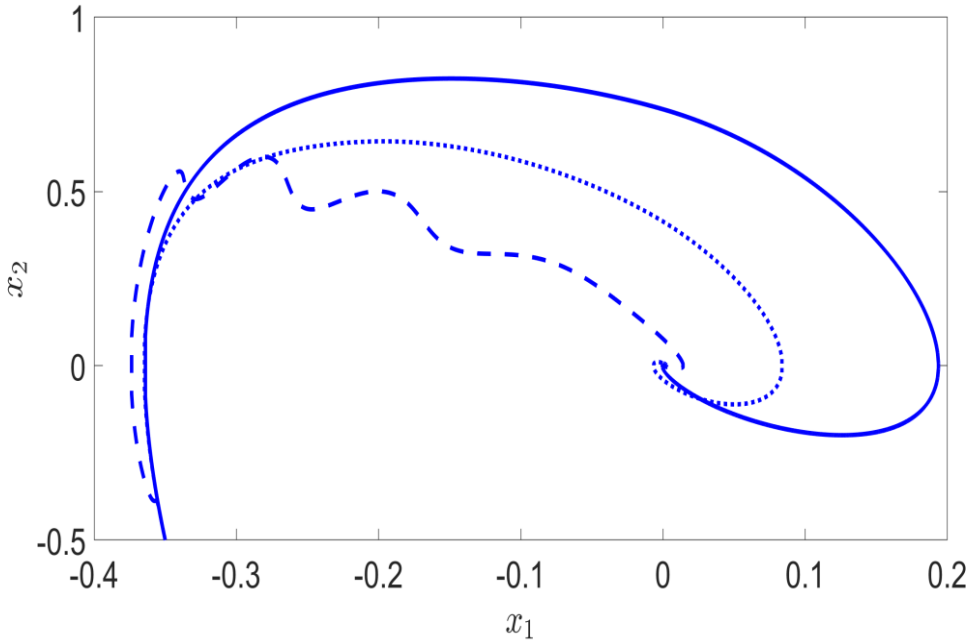

**Fig. 1:** The projection on the $x_1 - x_2$ plane of the solutions of the nominal closed-loop system (5.1) with (5.2) (solid line) with known $\theta$, the closed-loop system (5.1) with the classical extended-matching adaptive controller (5.12) (dashed line) and the closed-loop system (5.1) with the proposed event-triggered adaptive scheme (3.1), (3.2), (3.4), (3.5), (3.16) (dotted line) for the 1st initial condition. No measurement errors are present.

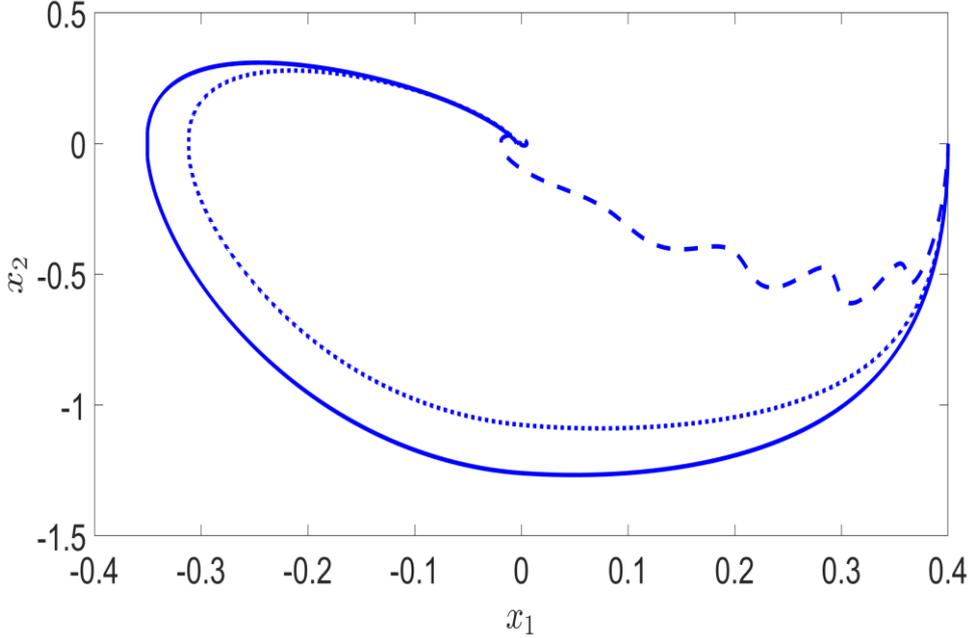

**Fig. 2:** The projection on the $x_1 - x_2$ plane of the solutions of the nominal closed-loop system (5.1) with (5.2) (solid line) with known $\theta$, the closed-loop system (5.1) with the classical extended-matching adaptive controller (5.12) (dashed line) and the closed-loop system (5.1) with the proposed event-triggered adaptive scheme (3.1), (3.2), (3.4), (3.5), (3.16) (dotted line) for the 2$^{nd}$ initial condition. No measurement errors are present.

The design based on the extended matching presents the smallest overshoot in the $x_1, x_2$ state components for both initial conditions as shown in Figures 1, 2 when no measurement errors are present. However, this is not the case for the entire state vector. Figures 3, 4 present the evolution of the norm of the state vector $|x(t)|$ for $t \in [0,4]$ and for both initial conditions, when no measurement errors are present. It is clear that the closed-loop system with the extended matching design exhibits very sharp overshoots for an initial transient period. Figures 5, 6 show the reason that explains this phenomenon (for the 2$^{nd}$ initial condition but this holds for both initial conditions): the spikes occur at the times when $\hat{\theta}_5(t)$ takes values close to the lowest allowable limit ($\kappa - \varepsilon = 0.999$). During this transient period the value of $\hat{\theta}_4(t)$ changes and assumes a stabilizing value (around 1.5), which allows the termination of the transient period.

We also studied the effect of measurement errors, i.e., the case where we measure

$$\hat{x}(t) = x(t) + e(t) \qquad (5.14)$$

where $e(t) \in \Re^3$ is the measurement error. In this case the nominal feedback becomes $u(t) = k(\theta, \hat{x}(t))$, the classical extended-matching adaptive controller is given by (5.12) with $x$ replaced by $\hat{x}$ and the proposed event-triggered adaptive scheme is given by (3.1), (3.2), (3.4), (3.5), (3.16), (5.7), (5.8), (5.9), (5.10) with $x$ replaced by $\hat{x}$. We used

$$e(t) = 0.01\sin(14\pi t)(1,1,1)', \text{ for } t \geq 0 \qquad (5.15)$$

and the results are shown in Fig. 9, Fig. 10, Fig. 11 and Fig. 12. It is clearly seen by Fig. 11 that the parameter estimation process by the BaLSI presents robustness with respect to measurement errors: the identifier manages to bring the parameter estimates very close to the exact parameter values in a few triggers even in the presence of measurement errors. However, Fig.9 and Fig.10 show that the overshoot exhibited by the proposed event-triggered adaptive scheme is given by (3.1), (3.2), (3.4), (3.5), (3.16), (5.7), (5.8), (5.9), (5.10) is larger when measurement errors are present due to the delayed convergence of the parameter estimates. Fig. 9, Fig. 10 and Fig. 12 show that the behavior of the closed-loop system with the classical extended-matching adaptive controller in the presence of measurement errors does not differ significantly from the behavior in the absence of measurement errors.

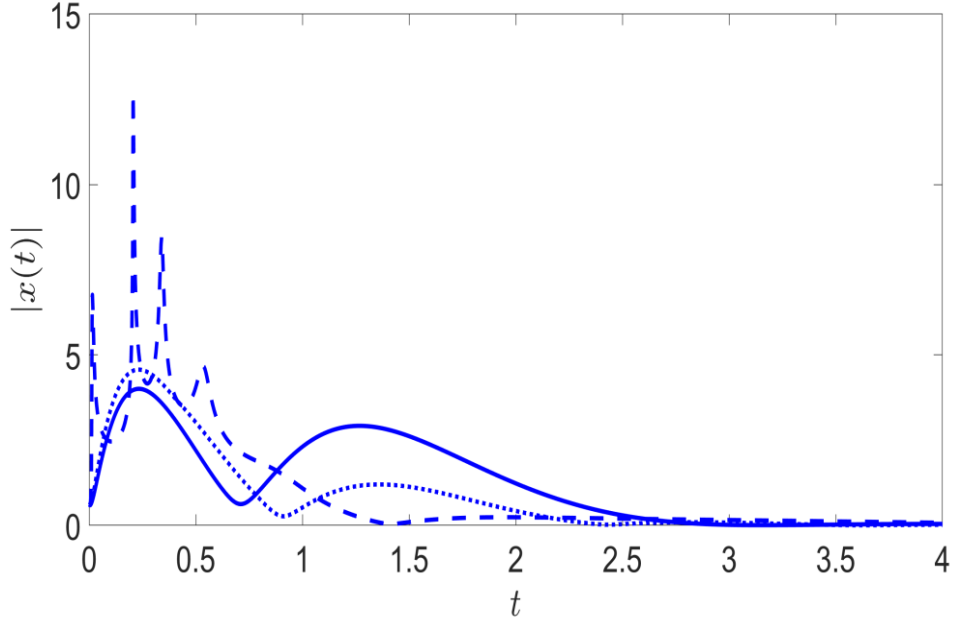

**Fig. 3:** The evolution of the value of $|x(t)|$ for $t \in [0,4]$ and 1st initial condition: solid line for the nominal closed-loop system (5.1) with (5.2) with known $\theta$, dashed line for the closed-loop system (5.1) with the classical extended-matching adaptive controller (5.12) and dotted line for the closed-loop system (5.1) with the proposed event-triggered adaptive scheme (3.1), (3.2), (3.4), (3.5), (3.16). No measurement errors are present.

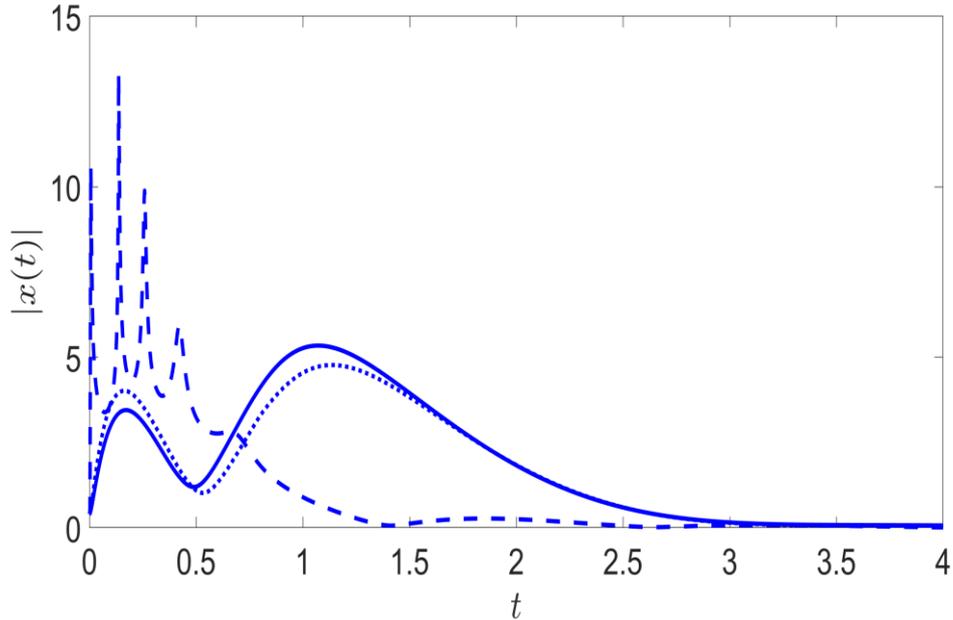

**Fig. 4:** The evolution of $|x(t)|$ for $t \in [0,4]$ and 2nd initial condition: solid line for the nominal closed-loop system (5.1) with (5.2) with known $\theta$, dashed line for the closed-loop system (5.1) with the classical extended-matching adaptive controller (5.12) and dotted line for the closed-loop system (5.1) with the proposed event-triggered adaptive scheme (3.1), (3.2), (3.4), (3.5), (3.16). No measurement errors are present.

The numerical results allow us to conclude that the proposed event-triggered adaptive scheme (3.1), (3.2), (3.4), (3.5), (3.15) exhibits exponential convergence of the state to zero and exact finite-time estimation of the unknown parameters (at least for all tested initial conditions and when measurement errors are absent) as well as robustness with respect to small measurement errors. However, these features come at a cost: the computational effort and the memory requirements for the implementation of the proposed event-triggered adaptive scheme (3.1), (3.2), (3.4), (3.5), (3.16) are significantly larger than those of the extended-matching design. Notice that the implementation of the proposed event-triggered adaptive scheme (3.1), (3.2), (3.4), (3.5), (3.16) requires 46 additional first-order ODEs to be solved in parallel to the three ODEs of the system, while the extended matching design requires only five additional first-order ODEs.

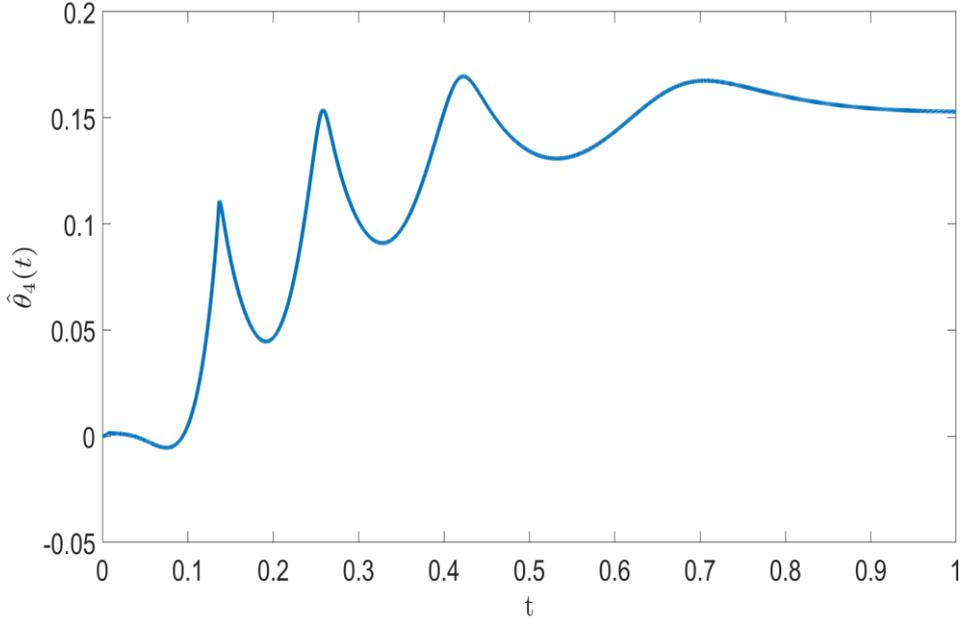

**Fig. 5:** The evolution of the value of $\hat{\theta}_4(t)$ for $t \in [0,1]$ and 2$^{nd}$ initial condition for the closed-loop system (5.1) with the classical extended-matching adaptive controller (5.12).

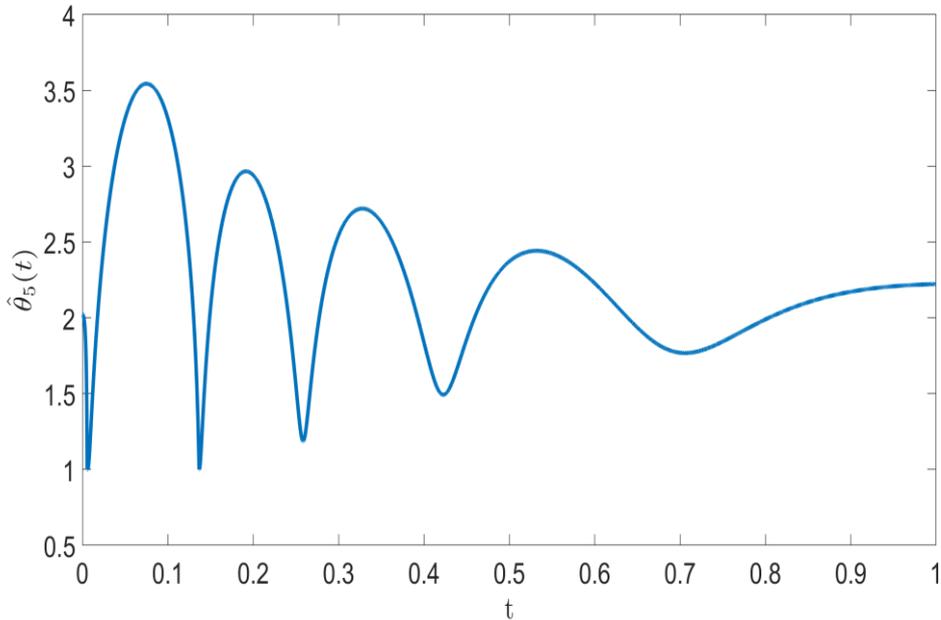

**Fig. 6:** The evolution of $\hat{\theta}_5(t)$ for $t \in [0,1]$ and 2$^{nd}$ initial condition for the closed-loop system (5.1) with the classical extended-matching adaptive controller (5.12).

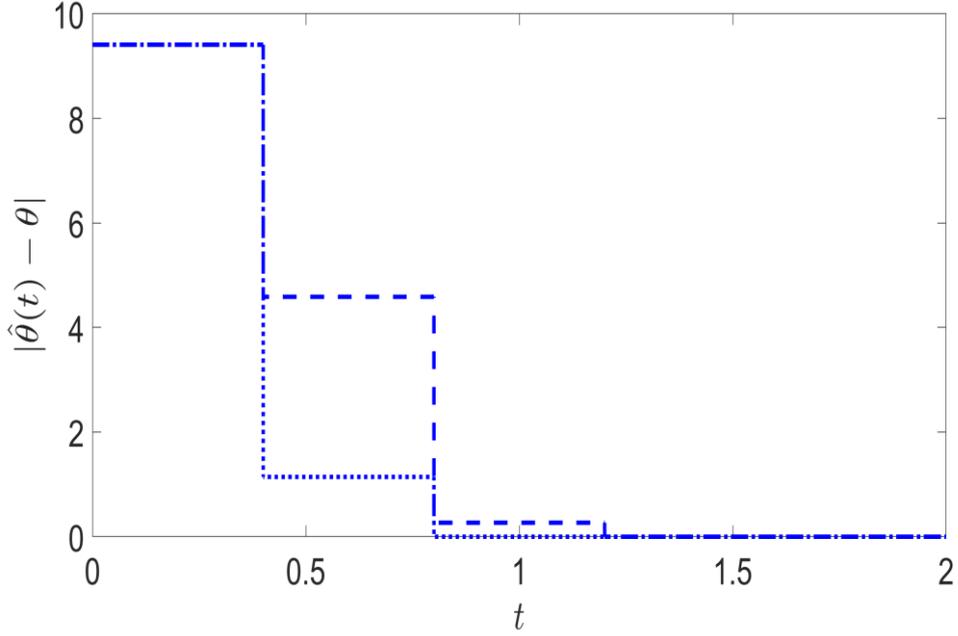

**Fig. 7:** The evolution of the Euclidean norm of the parameter estimation error $|\hat{\theta}(t) - \theta|$ for the closed-loop system (5.1) with the proposed event-triggered adaptive scheme (3.1), (3.2), (3.4), (3.5), (3.16): dashed line for the 1st initial condition and dotted line for the 2nd initial condition. No measurement errors are present.

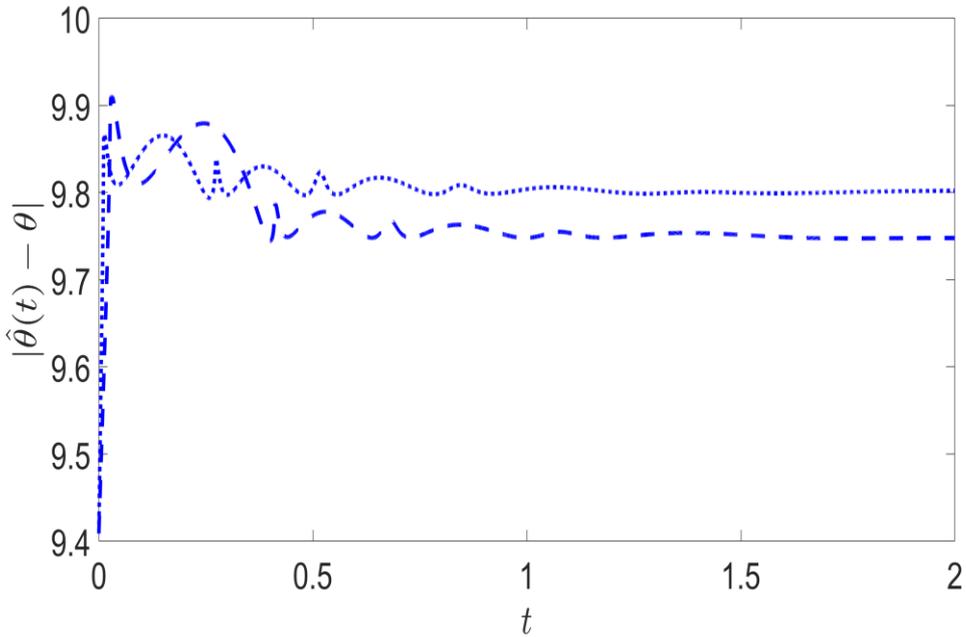

**Fig. 8:** The evolution of the Euclidean norm of the parameter estimation error $|\hat{\theta}(t) - \theta|$ for the closed-loop system (5.1) with the classical extended-matching adaptive controller (5.12): dashed line for the 1st initial condition and dotted line for the 2nd initial condition. No measurement errors are present.

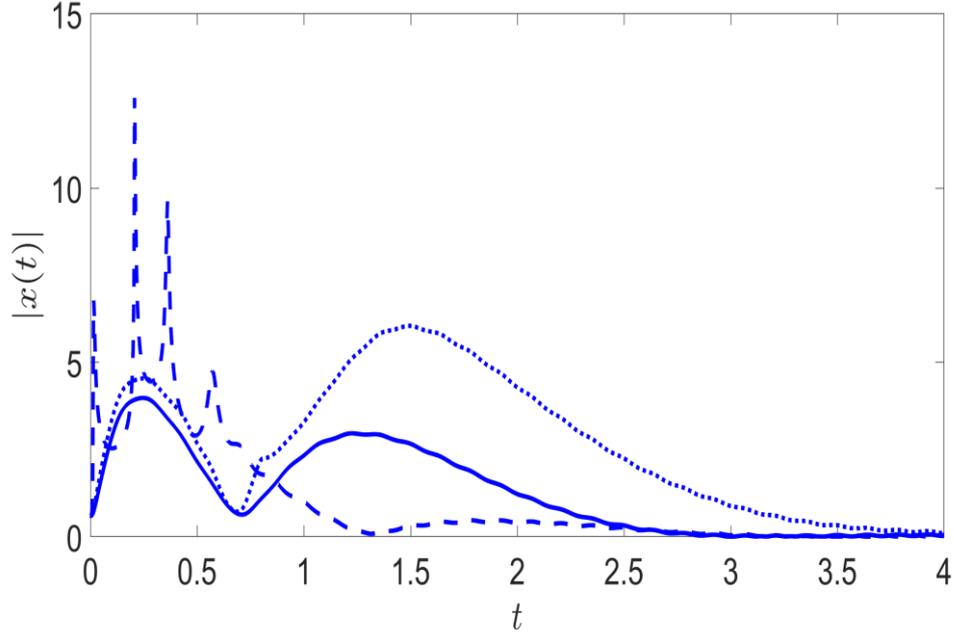

**Fig. 9:** The evolution of the value of $|x(t)|$ for $t \in [0,4]$ and 1$^{st}$ initial condition: solid line for the nominal closed-loop system (5.1) with (5.2) with known $\theta$, dashed line for the closed-loop system (5.1) with the classical extended-matching adaptive controller (5.12) and dotted line for the closed-loop system (5.1) with the proposed event-triggered adaptive scheme (3.1), (3.2), (3.4), (3.5), (3.16). Measurement errors given by (5.15).

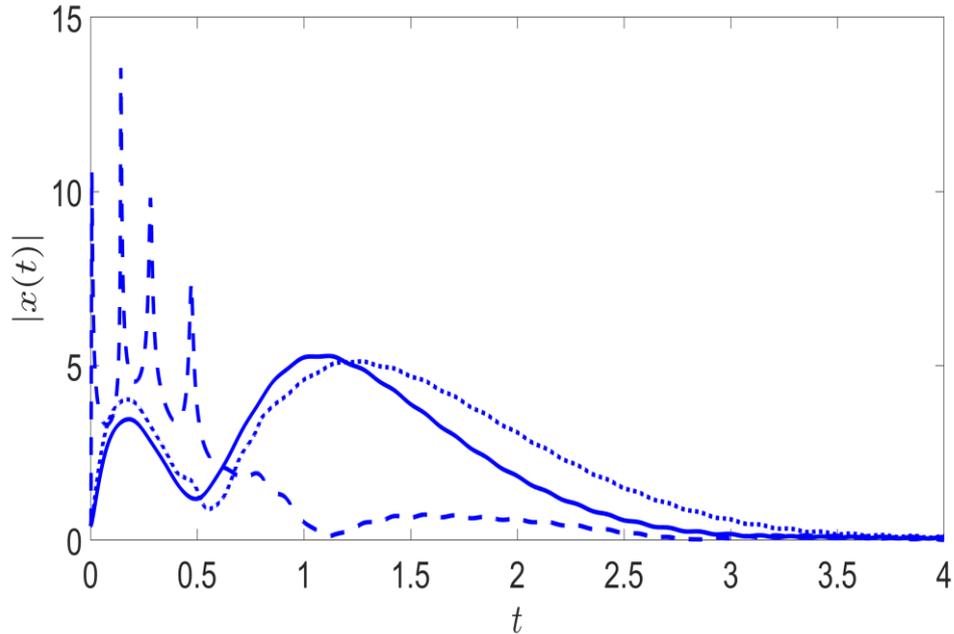

**Fig. 10:** The evolution of $|x(t)|$ for $t \in [0,4]$ and 2$^{nd}$ initial condition: solid line for the nominal closed-loop system (5.1) with (5.2) with known $\theta$, dashed line for the closed-loop system (5.1) with the classical extended-matching adaptive controller (5.12) and dotted line for the closed-loop system (5.1) with the proposed event-triggered adaptive scheme (3.1), (3.2), (3.4), (3.5), (3.16). Measurement errors given by (5.15).

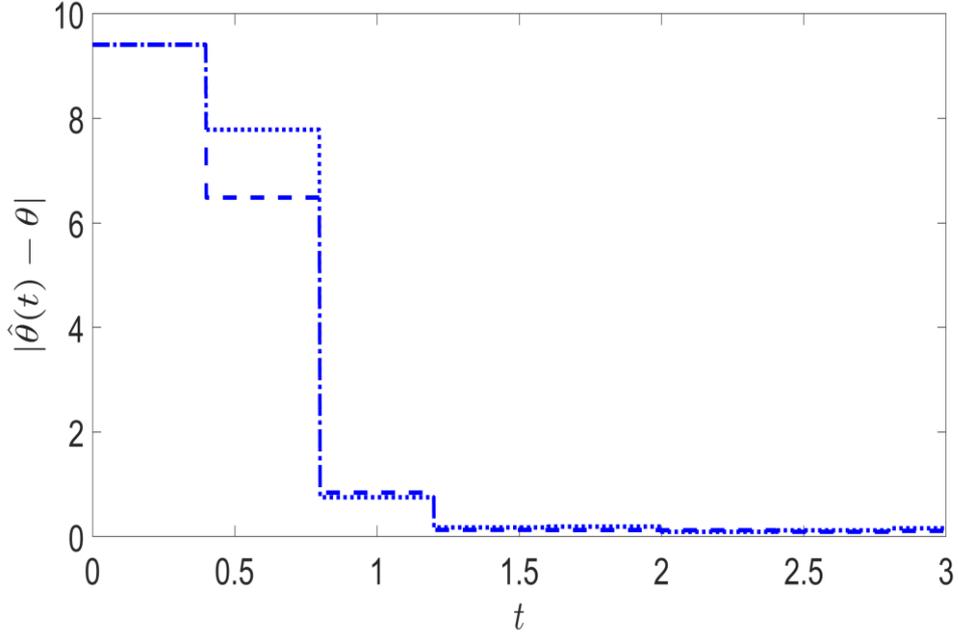

**Fig. 11:** The evolution of the Euclidean norm of the parameter estimation error $|\hat{\theta}(t) - \theta|$ for the closed-loop system (5.1) with the proposed event-triggered adaptive scheme (3.1), (3.2), (3.4), (3.5), (3.16): dashed line for the 1st initial condition and dotted line for the 2nd initial condition. Measurement errors given by (5.15).

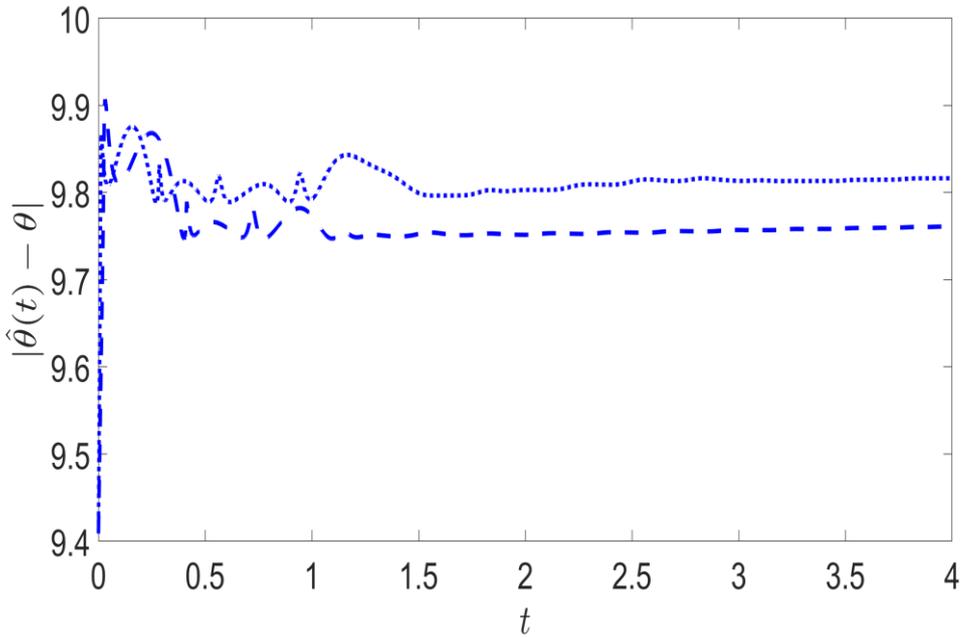

**Fig. 12:** The evolution of the Euclidean norm of the parameter estimation error $|\hat{\theta}(t) - \theta|$ for the closed-loop system (5.1) with the classical extended-matching adaptive controller (5.12): dashed line for the 1st initial condition and dotted line for the 2nd initial condition. Measurement errors given by (5.15).

# 6. Proofs of Results

We start this section with the proof of Theorem 4.1.

**Proof of Theorem 4.1:** The first claim is a direct consequence of the event trigger given by (3.4) and (3.5). The proof of the first claim is straightforward and is omitted.

**Claim 1:** *If a solution $(x(t),\hat{\theta}(t)) \in \Re^n \times \Theta$ of the closed-loop system (1.1) with (3.1), (3.2), (3.4), (3.5) and (3.13) is defined on $t \in [0,\tau_i]$ for certain $i \in Z_+$, then the solution is defined on $t \in [0,\tau_{i+1}]$. Moreover, it holds that*

$$V(\hat{\theta}(\tau_i), x(t)) \leq V(\hat{\theta}(\tau_i), x(\tau_i)) + a(x(\tau_i)), \text{ for all } t \in [\tau_i, \tau_{i+1}] \tag{6.1}$$

It follows from Claim 1 that for every initial condition, the corresponding solution $(x(t),\hat{\theta}(t)) \in \Re^n \times \Theta$ of the closed-loop system (1.1) with (3.1), (3.2), (3.4), (3.5) and (3.13) is defined for all $t \in \left[0, \sup_{i \geq 0}(\tau_i)\right)$. Uniqueness is a straightforward consequence of the assumed regularity properties for $f, g, k$.

The second claim clarifies what happens when the parameter estimation error becomes zero at the time of an event.

**Claim 2:** *If a solution $(x(t),\hat{\theta}(t)) \in \Re^n \times \Theta$ of the closed-loop system (1.1) with (3.1), (3.2), (3.4), (3.5) and (3.13) satisfies $\hat{\theta}(\tau_i) = \theta$ for certain $i \in Z_+$, then the solution is defined for all $t \geq 0$ and satisfies $\hat{\theta}(t) = \theta$ for all $t \geq \tau_i$ and $\tau_j = \tau_i + (j-i)T$ for all $j \geq i$.*

**Proof of Claim 2:** Notice that $\hat{\theta}(t) = \theta$ for all $t \in [\tau_i, \tau_{i+1})$. Assume first that $x(\tau_i) \neq 0$. Since for every $\theta \in \Theta$, $y_0 \in \Re^n$ the solution $y(t) \in \Re^n$ of $\dot{y} = f(y, k(\theta, y)) + g(y, k(\theta, y))\theta$ with initial condition $y(0) = y_0$ satisfies the inequality $V(\theta, y(t)) \leq V(\theta, y_0)$ for all $t \geq 0$ (recall assumption (H1)), it follows that $V(\hat{\theta}(\tau_i), x(t)) \leq V(\hat{\theta}(\tau_i), x(\tau_i)) < V(\hat{\theta}(\tau_i), x(\tau_i)) + a(x(\tau_i))$ for all $t \in [\tau_i, \tau_{i+1}]$. Consequently, it follows from (3.4) that $\tau_{i+1} = \tau_i + T$. The same conclusion follows from (3.5) if $x(\tau_i) = 0$. Since equation (3.10) holds, it follows from (3.13) that $\hat{\theta}(\tau_{i+1}) = \theta$. Applying the argument inductively, we conclude that $\hat{\theta}(t) = \theta$ for all $t \geq \tau_i$ and $\tau_j = \tau_i + (j-i)T$ for all $j \geq i$. The proof of Claim 2 is complete. ◁

We next define the following matrix for all $\tau \in \left[0, \sup_{i \geq 0}(\tau_i)\right)$:

$$G(\tau) := \left( \int_0^\tau \int_0^\tau q'(t,s) q(t,s) \, ds \, dt \right) \in \Re^{l \times l} \tag{6.2}$$

where $q(t,s)$ is defined by (3.9). For this matrix we can prove the following claim. Notice that $N(G(\tau))$ denotes the null space of $G(\tau)$, i.e.,

$$N(G(\tau)) := \left\{ \xi \in \Re^l : G(\tau)\xi = 0 \right\}.$$

**Claim 3:** *A vector $\xi \in \Re^l$ belongs to $N(G(\tau))$ for some $\tau \in \left(0, \sup_{i \geq 0}(\tau_i)\right]$ if and only if $q(t,s)\xi = 0$ for all $t,s \in [0,\tau]$ and $g(x(t),u(t))\xi = 0$ for all $t \in [0,\tau)$.*

**Proof of Claim 3:** If $G(\tau)\xi = 0$ then definition (6.2) implies that $\xi' G(\tau)\xi = \int_0^\tau \int_0^\tau |q(t,s)\xi|^2 \, ds \, dt = 0$. By continuity of the mapping $q(t,s)\xi$ for all $t,s \in [0,\tau]$, we obtain $q(t,s)\xi = 0$ for all $t,s \in [0,\tau]$. Using definition (3.9) and noticing that $t \to g(x(t),u(t))\xi$ is continuous for all $t \in [0,\tau] \setminus \bigcup_{i \geq 0} \{\tau_i\}$ and right-continuous for all $t \in [0,\tau)$, we obtain $g(x(t),u(t))\xi = 0$ for all $t \in [0,\tau)$. On the other hand if $q(t,s)\xi = 0$ for all $t,s \in [0,\tau]$ then definition (6.2) implies that $G(\tau)\xi = 0$. The proof of Claim 3 is complete. ◁

Claim 3 and the fact that $N(G(0)) = \Re^l$ (recall definition (6.2)) implies that the following inclusion holds for all $i \geq 0$:
$$N(G(\tau_{i+1})) \subseteq N(G(\tau_i)) \tag{6.3}$$

The following claim clarifies what happens when a switching in the value of the parameter estimate $\hat{\theta}(t)$ occurs.

**Claim 4:** *If $\hat{\theta}(\tau_{i+1}) \neq \hat{\theta}(\tau_i)$ then $\dim(N(G(\tau_{i+1}))) < \dim(N(G(\tau_i)))$.*

**Proof of Claim 4:** If $\hat{\theta}(\tau_{i+1}) \neq \hat{\theta}(\tau_i)$ then it follows from (3.13) that $G(\tau_{i+1})\hat{\theta}(\tau_i) \neq Z(\tau_{i+1})$, where
$$Z(\tau_{i+1}) := \int_0^{\tau_{i+1}} \int_0^{\tau_{i+1}} q'(t,s) p(t,s) \, ds \, dt \tag{6.4}$$

It follows from (3.10) and definitions (6.2), (3.13), that the vector $\xi = \hat{\theta}(\tau_i) - \theta$ satisfies $G(\tau_{i+1})\xi \neq 0$. Moreover, it follows from (3.10) and definitions (6.2), (3.13), that $G(\tau_i)\xi = 0$. Therefore, there exists a vector $\xi \in N(G(\tau_i))$ with $\xi \notin N(G(\tau_{i+1}))$.

We next prove by contradiction that $\dim(N(G(\tau_{i+1}))) < \dim(N(G(\tau_i)))$. Suppose that $\nu = \dim(N(G(\tau_{i+1}))) \geq \mu = \dim(N(G(\tau_i))) > 0$ (the case $\mu = 0$ can be excluded by Claim 2 and (3.10), (3.13) which would imply that $\hat{\theta}(\tau_{i+1}) = \hat{\theta}(\tau_i) = \theta$). Let $\{\xi_1,...,\xi_\mu\}$ be a basis for $N(G(\tau_i))$ and let $\{w_1,...,w_\nu\}$ be a basis for $N(G(\tau_{i+1}))$. At least one of the vectors $\xi_1,...,\xi_\mu$ must not belong to $N(G(\tau_{i+1}))$ because otherwise we would have $N(G(\tau_i)) \subseteq N(G(\tau_{i+1}))$, which contradicts the existence of a vector $\xi \in N(G(\tau_i))$ with $\xi \notin N(G(\tau_{i+1}))$. Let $j \in \{1,...,\mu\}$ with $\xi_j \notin N(G(\tau_{i+1}))$. Notice that the set of $\nu + 1 \geq \mu + 1$ vectors $\{w_1,...,w_\nu,\xi_j\}$ is a linearly independent set of vectors in $N(G(\tau_i))$, contradicting the assumption that $\mu = \dim(N(G(\tau_i)))$. The proof of Claim 4 is complete. ◁

A direct consequence of Claim 4 is the fact that at most $l$ switchings of the value of the parameter estimate $\hat{\theta}(t)$ can occur in the time interval $\left[0, \sup_{i \geq 0}(\tau_i)\right)$.

We next notice that (3.10) implies that the parameter update law (3.13) satisfies the following estimate for all $i \in Z_+$:
$$\left|\hat{\theta}(\tau_{i+1}) - \hat{\theta}(\tau_i)\right| \leq \left|\theta - \hat{\theta}(\tau_i)\right| \tag{6.5}$$

Using the triangle inequality and (6.5) we obtain for all $i \in Z_+$:

$$\left|\theta - \hat{\theta}(\tau_{i+1})\right| \leq 2\left|\theta - \hat{\theta}(\tau_i)\right| \tag{6.6}$$

Since at most $l$ switchings of the value of the parameter estimate $\hat{\theta}(t)$ can occur, we obtain from (6.6):

$$\left|\hat{\theta}(t) - \theta\right| \leq 2^l \left|\hat{\theta}(0) - \theta\right|, \text{ for all } t \in \left[0, \sup_{i \geq 0}(\tau_i)\right) \tag{6.7}$$

Let $s \geq 0$, $\theta \in \Theta$ be given and define the functions:

$$\begin{aligned} \bar{V}(x;\theta,s) &:= \min\left\{V(\vartheta,x): \vartheta \in \Theta, |\vartheta - \theta| \leq 2^l s\right\} \\ W(x;\theta,s) &:= \exp(2\sigma T)\left(\max\left\{V(\vartheta,x): \vartheta \in \Theta, |\vartheta - \theta| \leq 2^l s\right\} + a(x)\right) \end{aligned}, \text{ for } x \in \Re^n \tag{6.8}$$

Proposition 2.9 on page 21 in [6] implies that the mappings $x \to \bar{V}(x;\theta,s)$, $x \to W(x;\theta,s)$ are continuous and positive definite for each fixed $s \geq 0$, $\theta \in \Theta$. Moreover, assumption (H2) guarantees that for each fixed $s \geq 0$, $\theta \in \Theta$ the mapping $x \to \bar{V}(x;\theta,s)$ is radially unbounded. Consequently, Proposition 2.2 on page 107 in [19] implies that for each fixed $s \geq 0$, $\theta \in \Re^l$, there exist functions $a_{\theta,s} \in K_\infty$, $\beta_{\theta,s} \in K_\infty$ such that

$$a_{\theta,s}(|x|) \leq \bar{V}(x;\theta,s), \ \beta_{\theta,s}(|x|) \geq W(x;\theta,s), \text{ for all } x \in \Re^n \tag{6.9}$$

For any given $s \geq 0$, $\theta \in \Theta$, define the *KL* function:

$$\bar{\omega}_{\theta,s}(r,t) := a_{\theta,s}^{-1}\left(\exp(-2\sigma t)\beta_{\theta,s}(r,t)\right), \text{ for } t,r \geq 0 \tag{6.10}$$

The following claim clarifies what happens when $\hat{\theta}(\tau_{i+1}) = \hat{\theta}(\tau_i)$.

**Claim 5:** *If $\hat{\theta}(\tau_{i+1}) = \hat{\theta}(\tau_i)$ then*

$$\nabla V(\hat{\theta}(\tau_i),x(t))\dot{x}(t) \leq -2\sigma V(\hat{\theta}(\tau_i),x(t)), \text{ for all } t \in [\tau_i,\tau_{i+1}) \tag{6.11}$$

*Moreover, $V(\hat{\theta}(\tau_i),x(t)) \leq \exp(-2\sigma(t-\tau_i))V(\hat{\theta}(\tau_i),x(\tau_i))$ for all $t \in [\tau_i,\tau_{i+1}]$ and $\tau_{i+1} = \tau_i + T$.*

**Proof of Claim 5:** It follows from (3.10) and definitions (6.2), (3.13), that there exists a vector $\xi \in N(G(\tau_{i+1}))$ such that $\hat{\theta}(\tau_i) = \theta + \xi$. Moreover, it follows from Claim 3 that $g(x(t),u(t))\xi = 0$ for all $t \in [\tau_i,\tau_{i+1})$. Therefore, we obtain from (1.1) and (3.1):

$$\dot{x}(t) = f\left(x(t),k(\hat{\theta}(\tau_i),x(t))\right) + g\left(x(t),k(\hat{\theta}(\tau_i),x(t))\right)\hat{\theta}(\tau_i), \text{ for all } t \in [\tau_i,\tau_{i+1}) \tag{6.12}$$

Inequality (6.11) is a direct consequence of inequality (2.3) and equation (6.12). The fact that $V(\hat{\theta}(\tau_i),x(t)) \leq \exp(-2\sigma(t-\tau_i))V(\hat{\theta}(\tau_i),x(\tau_i))$ for all $t \in [\tau_i,\tau_{i+1}]$ follows from the differential inequality (6.11) and the fact that the mapping $t \to V(\hat{\theta}(\tau_i),x(t))$ is continuous. Finally, (3.4) implies that the event-trigger is not activated and thus (3.2), (3.5) give $\tau_{i+1} = \tau_i + T$. The proof of Claim 5 is complete. ◁

Using Claim 1, the fact that $\tau_{i+1} \leq \tau_i + T$ for all $i \geq 0$ (recall (3.2)) and Claim 5, we conclude that the following statements hold:

- If $\hat{\theta}(\tau_{i+j}) = \hat{\theta}(\tau_i)$ for $j=1,\ldots,\mu$, $\mu \geq 1$ then $V(\hat{\theta}(\tau_i), x(t)) \leq \exp(-2\sigma(t-\tau_i))V(\hat{\theta}(\tau_i), x(\tau_i))$ for $t \in [\tau_i, \tau_{i+\mu}]$ and $\tau_{i+\mu} = \tau_i + \mu T$.
- If $\hat{\theta}(\tau_{i+1}) \neq \hat{\theta}(\tau_i)$ then $V(\hat{\theta}(\tau_i), x(t)) \leq \exp(-2\sigma(t-\tau_i))\exp(2\sigma T)\left(V(\hat{\theta}(\tau_i), x(\tau_i)) + a(x(\tau_i))\right)$ for $t \in [\tau_i, \tau_{i+1}]$.

Using estimate (6.7) and definitions (6.8), we are in a position to rephrase the above statements as follows:

- If $\hat{\theta}(\tau_{i+j}) = \hat{\theta}(\tau_i)$ for $j=1,\ldots,\mu$, $\mu \geq 1$ then $\bar{V}(x(t); \theta, |\hat{\theta}(0) - \theta|) \leq \exp(-2\sigma(t-\tau_i))W(x(\tau_i); \theta, |\hat{\theta}(0) - \theta|)$ for $t \in [\tau_i, \tau_{i+\mu}]$ and $\tau_{i+\mu} = \tau_i + \mu T$.
- If $\hat{\theta}(\tau_{i+1}) \neq \hat{\theta}(\tau_i)$ then $\bar{V}(x(t); \theta, |\hat{\theta}(0) - \theta|) \leq \exp(-2\sigma(t-\tau_i))W(x(\tau_i); \theta, |\hat{\theta}(0) - \theta|)$ for $t \in [\tau_i, \tau_{i+1}]$.

Finally, using (6.9), the above statements and definition (6.10), we are in a position to guarantee the following facts:

**(F1)** If $\hat{\theta}(\tau_{i+j}) = \hat{\theta}(\tau_i)$ for $j=1,\ldots,\mu$, $\mu \geq 1$ then $|x(t)| \leq \bar{\omega}_{\theta,s}(|x(\tau_i)|, t-\tau_i)$ for $t \in [\tau_i, \tau_{i+\mu}]$ with $s = |\hat{\theta}(0) - \theta|$ and $\tau_{i+\mu} = \tau_i + \mu T$.

**(F2)** If $\hat{\theta}(\tau_{i+1}) \neq \hat{\theta}(\tau_i)$ then $|x(t)| \leq \bar{\omega}_{\theta,s}(|x(\tau_i)|, t-\tau_i)$ for $t \in [\tau_i, \tau_{i+1}]$ with $s = |\hat{\theta}(0) - \theta|$.

Since at most $l$ switchings of the value of the parameter estimate $\hat{\theta}(t)$ can occur, Fact (F1) implies that $\tau_i \geq (i-l)T$ for all $i \geq l$. Therefore, we conclude that $\sup_{i \geq 0}(\tau_i) = +\infty$ and that every solution $(x(t), \hat{\theta}(t)) \in \Re^n \times \Theta$ of the closed-loop system (1.1) with (3.1), (3.2), (3.4), (3.5) and (3.13) is defined for all $t \geq 0$.

Let $H(x(0), \hat{\theta}(0), \theta) \subseteq \{\tau_1, \tau_2, \ldots\}$ be the set of all times $\tau_i$ with $\hat{\theta}(\tau_i) \neq \hat{\theta}(\tau_{i-1})$. Notice that the set $H(x(0), \hat{\theta}(0), \theta)$ may be empty and can have at most $l$ members (since at most $l$ switchings of the value of the parameter estimate $\hat{\theta}(t)$ can occur). By virtue of Fact VI in [18], for any given $s \geq 0$, $\theta \in \Theta$, we are in a position to find a function $\Omega_{\theta,s} \in KL$, $i = 0,\ldots,l$, such that:

$$\sup_{0 \leq \tau \leq t}\left(\bar{\omega}_{\theta,s}(\bar{\omega}_{\theta,s}(r,\tau), t-\tau)\right) \leq \Omega_{\theta,s}(r,t), \text{ for all } t, r \geq 0 \text{ and } i = 0,\ldots,l-1 \quad (6.13)$$

Inequality (6.13) and Facts (F1), (F2) imply the following fact.

**(F3)** If $\tau_j \in H(x(0), \hat{\theta}(0), \theta)$ and $(\tau_i, \tau_j) \cap H(x(0), \hat{\theta}(0), \theta) = \varnothing$ then $|x(t)| \leq \Omega_{\theta,s}(|x(\tau_i)|, t-\tau_i)$ for $t \in [\tau_i, \tau_j]$ with $s = |\hat{\theta}(0) - \theta|$.

By virtue of Fact VI in [18], for any given $s \geq 0$, $\theta \in \Theta$, we are in a position to find functions $R_{i,\theta,s} \in KL$, $i = 0,\ldots,l$, such that:

$$\sup_{0 \leq \tau \leq t}\left(\Omega_{\theta,s}(R_{i,\theta,s}(r,\tau), t-\tau)\right) \leq R_{i+1,\theta,s}(r,t), \text{ for all } t, r \geq 0 \text{ and } i = 0,\ldots,l-1 \quad (6.14)$$

$$R_{0,\theta,s}(r,t) := \Omega_{\theta,s}(r,t), \text{ for all } t,r \geq 0 \tag{6.15}$$

We next show the following claim.

**Claim 6:** *If the cardinal number of the set $H(x(0),\hat{\theta}(0),\theta)$ is $\eta \in \{0,1,...,l\}$ then $|x(t)| \leq R_{\eta,\theta,s}(|x(0)|,t)$ for all $t \geq 0$ with $s = |\hat{\theta}(0) - \theta|$.*

**Proof of Claim 6:** By virtue of Fact (F1), the claim holds if $\eta = 0$. Indeed, if the cardinal number of the set $H(x(0),\hat{\theta}(0),\theta)$ is 0, i.e., if $\hat{\theta}(\tau_{i+1}) = \hat{\theta}(\tau_i)$ for all $i \geq 0$, then Fact (F1) implies that $|x(t)| \leq \bar{\omega}_{\theta,s}(|x(0)|,t)$ for all $t \geq 0$ with $s = |\hat{\theta}(0) - \theta|$. Definition (6.15) and the fact that $\bar{\omega}_{\theta,s}(r,t) \leq \Omega_{\theta,s}(r,t)$ for all $t,r \geq 0$ (a consequence of (6.13) and the fact that definitions (6.8), (6.10) guarantee that $r \leq \bar{\omega}_{\theta,s}(r,0)$) shows that the claim holds if $\eta = 0$.

Next suppose that $\eta > 0$. Let $T_1,...,T_\eta$ be such that $H(x(0),\hat{\theta}(0),\theta) = \{T_1,...,T_\eta\}$. By virtue of facts (F1) and (F3) we get:
- $|x(t)| \leq \Omega_{\theta,s}(|x(0)|,t)$ for $t \in [0,T_1]$ with $s = |\hat{\theta}(0) - \theta|$.
- $|x(t)| \leq \Omega_{\theta,s}(|x(T_i)|,t-T_i)$ for $t \in [T_i,T_{i+1}]$, $i = 1,...,\eta-1$ with $s = |\hat{\theta}(0) - \theta|$ (this case applies only when $\eta > 1$).
- $|x(t)| \leq \Omega_{\theta,s}(|x(T_\eta)|,t-T_\eta)$ for $t \geq T_\eta$ with $s = |\hat{\theta}(0) - \theta|$.

Combining all cases above with inequalities (6.14) and using the fact that $R_{i,\theta,s}(r,\tau) \leq R_{i+1,\theta,s}(r,t)$, for all $t,r \geq 0$ and $i = 0,...,l-1$, we get the desired estimate $|x(t)| \leq R_{\eta,\theta,s}(|x(0)|,t)$ for all $t \geq 0$ with $s = |\hat{\theta}(0) - \theta|$. The proof of Claim 6 is complete. ◁

Using the fact that the cardinal number of the set $H(x(0),\hat{\theta}(0),\theta)$ is at most $l$ and the fact that $R_{i,\theta,s}(r,\tau) \leq R_{i+1,\theta,s}(r,t)$, for all $t,r \geq 0$ and $i = 0,...,l-1$, we conclude that the required estimate $|x(t)| \leq \omega_{\theta,\hat{\theta}_0}(|x_0|,t)$ for all $t \geq 0$ holds for the solution of the hybrid closed-loop system (1.1), (3.1), (3.2), (3.4), (3.5), (3.13) and initial conditions $x(0) = x_0$, $\hat{\theta}(0) = \hat{\theta}_0$ with $\omega_{\theta,\hat{\theta}_0}(r,t) := R_{l,\theta,s}(r,t)$ and $s = |\hat{\theta}_0 - \theta|$.

Let $\tau > 0$ be the maximum time in the set $H(x(0),\hat{\theta}(0),\theta)$ when $H(x(0),\hat{\theta}(0),\theta) \neq \emptyset$ and $\tau = 0$ when $H(x(0),\hat{\theta}(0),\theta) = \emptyset$. It follows that $\hat{\theta}(t) = \theta_s = \hat{\theta}(\tau)$ for all $t \geq \tau$. If $\theta_s \neq \theta$ then it follows from (3.10) and definitions (6.2), (3.13), that $\xi = (\hat{\theta}(\tau) - \theta) \in N(G(\tau_{i+1}))$ for all $i \geq 0$ with $\tau_{i+1} > \tau$. It follows from Claim 3 that $g(x(t),u(t))\xi = 0$ for all $t \geq 0$.

Finally, if the parameter observability assumption (H3) holds then we can repeat all arguments in the proof of Theorem 4.1 in [20] and show that $\hat{\theta}(t) = \theta$ for all $t \geq NT$ when $x(0) \neq 0$.

The proof is complete. ◁

We next continue with the proof of Theorem 4.2.

**Proof of Theorem 4.2:** Since all assumptions of Theorem 4.1 hold for Theorem 4.2, the proof of Theorem 4.2 starts at the point where the proof of Theorem 4.1 ended. Therefore, all relations and everything written in the proof of Theorem 4.1 holds.

Define
$$A := \sup\{|x|^{-2} a(x) : x \neq 0, |x| \leq \delta\}. \tag{6.16}$$

Notice that due to (4.1), (6.16), (6.7) and definitions (6.8) for every $\theta \in \Theta$, $s \geq 0$ there exist constants $R^{\theta,s} > 0$, $K_2^{\theta,s} > K_1^{\theta,s} > 0$ such that

$$K_1^{\theta,s} |x|^2 \leq \bar{V}(x;\theta,s), \quad \exp(2\sigma T)(K_2^{\theta,s} + A)|x|^2 \geq W(x;\theta,s),$$
$$\text{for all } x \in \Re^n \text{ with } |x| \leq \min(R^{\theta,s}, \delta) \tag{6.17}$$

It follows from (6.17) and facts (F1), (F2) that the following statements hold.

**(S1)** If $\hat{\theta}(\tau_{i+j}) = \hat{\theta}(\tau_i)$ for $j = 1,\ldots,\mu$, $\mu \geq 1$ and $|x(\tau_i)| \leq \Gamma_{\theta,s}^{-1} \min(R^{\theta,s}, \delta)$, where $\Gamma_{\theta,s} := \exp(\sigma T)\sqrt{(K_1^{\theta,s})^{-1}(K_2^{\theta,s} + A)}$, then $|x(t)| \leq \Gamma_{\theta,s} \exp(-\sigma(t-\tau_i))|x(\tau_i)|$ for $t \in [\tau_i, \tau_{i+\mu}]$ with $s = |\hat{\theta}(0) - \theta|$.

**(S2)** If $\hat{\theta}(\tau_{i+1}) \neq \hat{\theta}(\tau_i)$ and $|x(\tau_i)| \leq \Gamma_{\theta,s}^{-1} \min(R^{\theta,s}, \delta)$, then $|x(t)| \leq \Gamma_{\theta,s} \exp(-\sigma(t-\tau_i))|x(\tau_i)|$ for $t \in [\tau_i, \tau_{i+1}]$ with $s = |\hat{\theta}(0) - \theta|$, where $\Gamma_{\theta,s} := \exp(\sigma T)\sqrt{(K_1^{\theta,s})^{-1}(K_2^{\theta,s} + A)}$.

Facts (S1), (S2) allow us to state the following fact.

**(S3)** If $\tau_j \in H(x(0), \hat{\theta}(0), \theta)$, $(\tau_i, \tau_j) \cap H(x(0), \hat{\theta}(0), \theta) = \emptyset$ and $|x(\tau_i)| \leq \Gamma_{\theta,s}^{-2} \min(R^{\theta,s}, \delta)$ then $|x(t)| \leq \Gamma_{\theta,s}^2 \exp(-\sigma(t-\tau_i))|x(\tau_i)|$ for $t \in [\tau_i, \tau_j]$ with $s = |\hat{\theta}(0) - \theta|$.

Using the fact that the cardinal number of the set $H(x(0), \hat{\theta}(0), \theta)$ is at most $l$ and facts (S1), (S3), we conclude that $|x(t)| \leq \Gamma_{\theta,s}^{2l} \exp(-\sigma t)|x(0)|$ for all $t \geq 0$, provided that $|x(0)| \leq \Gamma_{\theta,s}^{-2l} \min(R^{\theta,s}, \delta)$ with $\Gamma_{\theta,s} := \exp(\sigma T)\sqrt{(K_1^{\theta,s})^{-1}(K_2^{\theta,s} + A)}$ and $s = |\hat{\theta}_0 - \theta|$.

It follows that the desired estimate $|x(t)| \leq M_{\theta,\hat{\theta}_0} \exp(-\sigma t)|x_0|$ for all $t \geq 0$ and for all $\theta \in \Theta$, $x_0 \in \Re^n$, $\hat{\theta}_0 \in \Theta$ with $|x_0| \leq \bar{R}_{\theta,\hat{\theta}_0}$ holds for the solution of the hybrid closed-loop system (1.1) with (3.1), (3.2), (3.4), (3.5), (3.13) and initial conditions $x(0) = x_0$, $\hat{\theta}(0) = \hat{\theta}_0$ with $M_{\theta,\hat{\theta}_0} := \Gamma_{\theta,s}^{2l}$, $\bar{R}_{\theta,\hat{\theta}_0} := \Gamma_{\theta,s}^{-2l} \min(R^{\theta,s}, \delta)$ with $\Gamma_{\theta,s} := \exp(\sigma T)\sqrt{(K_1^{\theta,s})^{-1}(K_2^{\theta,s} + A)}$ and $s = |\hat{\theta}_0 - \theta|$.

The proof is complete. ◁

**Proof of Theorem 4.3:** The proof of Theorem 4.3 is almost identical with the proof of Theorem 4.2, except of the fact that no restrictions in the magnitude of $|x|$ are needed for the derivation of all estimates. ◁

**Proof of Corollary 4.4:** Let $\sigma > 0$ be given (arbitrary). The existence of continuous mappings $\Theta \ni \theta \to P(\theta) \in \Re^{n \times n}$, $\Theta \ni \theta \to \tilde{k}(\theta) \in \Re^n$ for which $P(\theta) = \{ p_{i,j}(\theta) : i,j = 1,\ldots,n \} \in \Re^{n \times n}$ is a symmetric, positive definite matrix and for which the inequality $\sum_{i=1}^{n}\sum_{j=1}^{n} p_{i,j}(\theta)\dot{x}_i x_j + \sum_{i=1}^{n}\sum_{j=1}^{n} p_{i,j}(\theta) x_i \dot{x}_j \leq -2\sigma \sum_{i=1}^{n}\sum_{j=1}^{n} p_{i,j}(\theta) x_i x_j$ with $u = \tilde{k}'(\theta) x$ holds for all $x \in \Re^n$ can be proved by induction on $n$. Indeed, for $n=1$ the statement holds with $p_{1,1}(\theta) = \frac{1}{2}$ and $\tilde{k}_1(\theta) = -\frac{\theta_{1,1} + \sigma}{\theta_{1,2}}$. If we assume that the statement holds for certain integer $n \geq 1$ then straightforward manipulations show that the statement is also true for $n+1$ with

$$p_{i,j}(\theta) = p_{i,j}(\theta) + \frac{1}{2}\tilde{k}_i(\theta)\tilde{k}_j(\theta) \quad i,j = 1,\ldots,n$$

$$p_{i,n+1}(\theta) = p_{n+1,i}(\theta) = -\frac{1}{2}\tilde{k}_i(\theta) \quad i = 1,\ldots,n$$

$$p_{n+1,n+1}(\theta) = \frac{1}{2}$$

$$\tilde{k}_1(\theta) = \theta_{n+1,n+2}^{-1}\left( \sum_{j=1}^{n} \tilde{k}_j(\theta)\theta_{j,1} - 2\theta_{n,n+1} p_{n,1}(\theta) - \theta_{n+1,1} \right) + \sigma \tilde{k}_1(\theta)$$

$$\tilde{k}_i(\theta) = \theta_{n+1,n+2}^{-1}\left( \sum_{j=i}^{n} \tilde{k}_j(\theta)\theta_{j,i} + \tilde{k}_{i-1}(\theta)\theta_{i-1,i} - 2\theta_{n,n+1} p_{n,i}(\theta) - \theta_{n+1,i} + \sigma \tilde{k}_i(\theta) \right) \quad i = 2,\ldots,n$$

$$\tilde{k}_{n+1}(\theta) = \theta_{n+1,n+2}^{-1}\left( \tilde{k}_n(\theta)\theta_{n,n+1} - \sigma - \theta_{n+1,n+1} \right)$$

The rest of the proof is a direct consequence of Theorem 4.3 with $f(x,u) \equiv 0$,

$$g(x,u)\theta := \begin{bmatrix} \theta_{1,1} x_1 + \theta_{1,2} x_2 \\ \vdots \\ \theta_{n,1} x_1 + \ldots + \theta_{n,n} x_n + \theta_{n,n+1} u \end{bmatrix}, \quad a(x) := A|x|^2, \quad V(\theta,x) := x'P(\theta)x \text{ and } k(\theta,x) := \tilde{k}'(\theta)x.$$ More specifically, inequality (4.2) is a consequence of continuity of the mapping $\Theta \ni \theta \to P(\theta) \in \Re^{n \times n}$, which implies that the function $V(\theta,x) = x'P(\theta)x$ is continuous on the compact set $\bar{\Theta} \times \{ x \in \Re^n : |x| = 1 \}$. The proof is complete. ◁

## 7. Concluding Remarks

The present work showed that regulation-triggered, adaptive schemes can guarantee exponential regulation even in the absence of persistency of excitation and parameter observability. The proposed adaptive scheme guarantees that the closed-loop system follows the trajectories of the nominal closed-loop system as well as a *KL* estimate for the state component $x$. However, a number of issues remain open:
  (i) the use of weighting functions in the BaLSI, and
  (ii) the numerical implementation of the BaLSI.
Both issues may be important in practice. An additional issue that should be addressed in future research is the study of sensitivity with respect to modeling and measurement errors. Although preliminary studies (see [20]) have shown important robustness properties with respect to various perturbations (vanishing and non-vanishing), the issue requires further (both theoretical and numerical) study.

**References**

[1] Anfinsen, H. and O. M. Aamo, "Adaptive Control of Linear 2×2 Hyperbolic Systems", *Automatica*, 87, 2018, 69-82.

[2] Arzen, K.-E., "A Simple Event-Based PID Controller", *Proceedings of the 1999 IFAC World Congress*, 1999, 423-428.

[3] Bernard, P. and M. Krstic, "Adaptive Output-feedback Stabilization of Non-Local Hyperbolic PDEs", *Automatica*, 50, 2014, 2692-2699.

[4] De Persis C., R. Saile and F. Wirth, "Parsimonious Event Triggered Distributed Control: A Zeno Free Approach", *Automatica*, 49, 2013, 2116-2124.

[5] Donkers T., and M. Heemels, "Output-Based Event-Triggered Control With Guaranteed L1-Gain and Improved and Decentralized Event-Triggering", *IEEE Transactions on Automatic Control*, 57, 2012, 1362-1376.

[6] Freeman, R. A. and P. V. Kokotovic, *Robust Nonlinear Control Design- State Space and Lyapunov Techniques*, Birkhauser, Boston, 1996.

[7] Garcia E., and P. J. Antsaklis, "Model-Based Event-Triggered Control for Systems with Quantization and Time-Varying Network Delays", *IEEE Transactions on Automatic Control*, 58, 2013, 422-434.

[8] Gawthrop P. J., and L. Wang, "Event-Driven Intermittent Control", *International Journal of Control*, 82, 2009, 2235-2248.

[9] Guinaldo, M., D. V. Dimarogonas, K. H. Johansson, J. Sanchez, and S. Dormido, "Distributed Event-Based Control Strategies for Interconnected Linear Systems", *IET Control Theory and Applications*, 7, 2013, 877-886.

[10] Heemels, W. P. M. H., J. H. Sandee, P. P. J. Van Den Bosch, "Analysis of Event-Driven Controllers for Linear Systems", *International Journal of Control*, 81, 2008, 571-590.

[11] Heemels, W. P. M. H., K. H. Johansson, and P. Tabuada, "An Introduction to Event-Triggered and Self-Triggered Control", *Proceedings of the 51st IEEE Conference in Decision and Control*, 2012, 3270-3285.

[12] Heemels, W. P. M. H., M. C. F. Donkers, and A. R. Teel, "Periodic Event-Triggered Control for Linear Systems", *IEEE Transactions on Automatic Control*, 58, 2013, 847-861.

[13] Hespanha, J. P., and A. S. Morse, "Certainty Equivalence Implies Detectability", *Systems and Control Letters*, 36, 1999, 1-13.

[14] Hespanha, J. P., D. Liberzon, and A. S. Morse, "Supervision of Integral-Input-to-State Stabilizing Controllers", *Automatica*, 38, 2002 1327-1335.

[15] Hespanha, J. P., D. Liberzon, and A. S. Morse, "Hysteresis-Based Switching Algorithms for Supervisory Control of Uncertain Systems", *Automatica*, 39, 2003, 263-272.

[16] Hespanha, J. P., D. Liberzon, and A. S. Morse, "Overcoming the Limitations of Adaptive Control by Means of Logic-Based Switching", *Systems and Control Letters*, 49, 2003, 49-65.

[17] Ioannou, P. A., and J. Sun, *Robust Adaptive Control*, Englewood Cliffs, NJ, USA: Prentice Hall, 1996.

[18] Karafyllis, I. and J. Tsinias, "Non-Uniform in Time Input-to-State Stability and the Small-Gain Theorem", *IEEE Transactions on Automatic Control*, 49, 2004, 196-216.

[19] Karafyllis, I., and Z.-P. Jiang, *Stability and Stabilization of Nonlinear Systems*, Springer-Verlag London (Series: Communications and Control Engineering), 2011.

[20] Karafyllis, I. and M. Krstic, "Adaptive Certainty-Equivalence Control With Regulation-Triggered Finite-Time Least-Squares Identification", *IEEE Transactions on Automatic Control*, 63, 2018, 3261-3275.

[21] Karafyllis, I., M. Krstic and K. Chrysafi, "Adaptive Boundary Control of Constant-Parameter Reaction-Diffusion PDEs Using Regulation-Triggered Finite-Time Identification", submitted to *Automatica* (see also arXiv:1805.08404 [math.OC]).

[22] Khalil, H.K., *Nonlinear Systems*, 2nd Edition, Prentice-Hall, 1996.

[23] Krstic, M., I. Kanellakopoulos and P. Kokotovic, "Adaptive Nonlinear Control Without Overparametrization", *Systems and Control Letters*, 19, 1992, 177-185.

[24] Krstic, M., I. Kanellakopoulos and P. Kokotovic, *Nonlinear and Adaptive Control Design*, Wiley, 1995.

[25] Liu, T., and Z. P. Jiang, "A Small-Gain Approach to Robust Event-Triggered Control of Nonlinear Systems", *IEEE Transactions on Automatic Control*, 2015, 60, 2015, 2072-2085.

[26] Lunze, J., and D. Lehmann, "A State-Feedback Approach to Event-Based Control", *Automatica*, 46, 2010, 211-215.

[27] Marchand, N, S. Durand, and J. F. G. Castellanos, "A General Formula for Event-Based Stabilization of Nonlinear Systems", *IEEE Transactions on Automatic Control*, 58, 2013, 1332-1337.

[28] Miller, D., E., and E. J. Davison, "An Adaptive Controller Which Provides Lyapunov Stability", *IEEE Transactions on Automatic Control*, 34, 1989, 599-609.

[29] Miller, D., E., and E. J. Davison, "An Adaptive Controller Which Provides an Arbitrarily Good Transient and Steady-State Response", *IEEE Transactions on Automatic Control*, 36, 1991, 68-81.

[30] Miller, D., E., "Near Optimal LQR Performance for a Compact Set of Plants", *IEEE Transactions on Automatic Control*, 51, 2006, 1423-1439.

[31] Miller, D., E., and J. R. Vale, "Pole Placement Adaptive Control With Persistent Jumps in the Plant Parameters", *Mathematics of Control, Signals, and Systems*, 26, 2014, 177–214.

[32] Monahemi, M. M., J. B. Barlow, and M. Krstic, "Control of Wing Rock Motion of Slender Delta Wings Using Adaptive Feedback Linearization", *Proceedings of AIAA Guidance, Navigation, and Control Conference*, Baltimore, MD, 1995, 1564-1575.

[33] Morse, A., S., "Supervisory Control of Families of Linear Set-Point Controllers-Part 1: Exact Matching", *IEEE Transactions on Automatic Control*, 41, 1996, 1413–1431.

[34] Morse, A., S., "Supervisory Control of Families of Linear Set-Point Controllers-Part 2: Robustness", *IEEE Transactions on Automatic Control*, 42, 1997, 1500–1515.

[35] Praly, L., G. Bastin, J.-B. Pomet, and Z. P. Jiang, "Adaptive Stabilization of Nonlinear Systems", *Foundations of Adaptive Control*, P. V. Kokotovic (Ed.), Springer-Verlag, 1991, 347-434.

[36] Sahoo, A., H. Xu and S. Jagannathan, "Neural Network-Based Adaptive Event-Triggered Control of Nonlinear Continuous-Time Systems", *Proceedings of the 2013 IEEE International Symposium on Intelligent Control (ISIC),* 2013, 35-40.

[37] Smyshlyaev, A. and M. Krstic, *Adaptive Control of Parabolic PDEs*, Princeton University Press, 2010.

[38] Sontag, E.D., *Mathematical Control Theory*, 2nd Edition, Springer-Verlag, New York, 1998.

[39] Tabuada, P., "Event-Triggered Real-Time Scheduling of Stabilizing Control Tasks", *IEEE Transactions on Automatic Control*, 52, 2007, 1680-1685.

[40] Tallapragada, P, and N. Chopra, "On Event Triggered Tracking for Nonlinear Systems", *IEEE Transactions on Automatic Control*, 58, 2013, 2343-2348.

[41] Tsinias, J., "Remarks on Asymptotic Controllability and Sampled-Data Feedback Stabilization for Autonomous Systems", *IEEE Transactions on Automatic Control*, 3, 2010, 721-726.

[42] Tsinias, J., "New Results on Sampled-Data Feedback Stabilization for Autonomous Nonlinear Systems", *Systems and Control Letters*, 61, 2012, 1032-1040.

[43] Vamvoudakis, K. G., "Event-Triggered Optimal Adaptive Control Algorithm for Continuous-Time Nonlinear Systems", *IEEE/CAA Journal of Automatica Sinica*, 1, 2014, 282-293.

[44] Wang, X., and N. Hovakimyan, "$L^1$ Adaptive Control of Event-Triggered Networked Systems", *Proceedings of the American Control Conference,* 2010, 2458-2463.

[45] Wang, X., and M. D. Lemmon, "Event-Triggering in Distributed Networked Control Systems", *IEEE Transactions on Automatic Control*, 56, 2011, 586-601.

[46] Xing, L., C. Wen, Z. Liu, H. Su and J. Cai, "Event-Triggered Adaptive Control for a Class of Uncertain Nonlinear Systems", *IEEE Transactions on Automatic Control*, 62, 2017, 2071-2076.

[47] Zhong, X. and H. He, "An Event-Triggered ADP Control Approach for Continuous-Time System With Unknown Internal States", *IEEE Transactions on Cybernetics*, 47, 2017, 683-694.

[48] Zhu, W., and Z. P. Jiang, "Event-Based Leader-Following Consensus of Multi-Agent Systems With Input Time Delay", *IEEE Transactions on Automatic Control*, 60, 2015, 1362-1367.

[49] Zhu, Y., D. Zhao, H. He and J. Ji, "Event-Triggered Optimal Control for Partially Unknown Constrained-Input Systems via Adaptive Dynamic Programming", *IEEE Transactions on Industrial Electronics*, 64, 2017, 4101-4109.